\documentclass{amsart}

\usepackage{graphicx}
\usepackage{latexsym}
\usepackage{amssymb, latexsym}
\usepackage{amsmath}
\usepackage{xcolor}
\usepackage{bbm}
\usepackage{mathrsfs}
\usepackage{amsxtra}
\usepackage{amsthm}
\usepackage{stmaryrd}
\usepackage{url}
\usepackage{verbatim}
\usepackage{xspace}
\usepackage{setspace}
\usepackage[hidelinks]{hyperref}
\usepackage{cleveref}

\newtheorem{theorem}{Theorem}[section]
\newtheorem{lemma}[theorem]{Lemma}
\newtheorem{proposition}[theorem]{Proposition}
\newtheorem{corollary}[theorem]{Corollary}

\theoremstyle{definition}
\newtheorem{definition}[theorem]{Definition}
\newtheorem{remark}[theorem]{Remark}
\newtheorem{question}[theorem]{Question}

\newcommand{\myref}{\Cref}

\newtheoremstyle{TheoremNum}
        {}{}           				   
        {\itshape}                      
        {}                              
        {\bfseries}                     
        {.}                             
        { }                             
        {\thmname{#1}\thmnote{ \bfseries #3}}
    \theoremstyle{TheoremNum}
    \newtheorem{duplicate}{Theorem}

\newcommand{\cminus}{\raisebox{.33\height}{\scalebox{0.75}{\ensuremath{-}}}}
\newcommand{\uminus}{\raisebox{.15\height}{\scalebox{0.75}{\ensuremath{-}}}}

\newcommand{\ZF}{{\rm ZF}}
\newcommand{\ZFC}{{\rm ZFC}}
\newcommand{\GB}{{\rm GB}}
\newcommand{\GBC}{{\rm GBC}}
\newcommand{\GBc}{{\rm GBc}}
\newcommand{\ZFminusrep}{\ensuremath{\ZF\cminus}\xspace}
\newcommand{\ZFCminusrep}{\ensuremath{\ZFC\cminus}\xspace}
\newcommand{\ZFCminus}{\ensuremath{\ZFC^{\uminus}}\xspace}
\newcommand{\ZFminus}{\ensuremath{\ZF^{\uminus}}\xspace}
\newcommand{\GBCminus}{\ensuremath{\GBC^{\uminus}}\xspace}
\newcommand{\GBminus}{\ensuremath{\GB^{\uminus}}\xspace}
\newcommand{\GBcminus}{\ensuremath{\GBc^{\uminus}}\xspace}

\newcommand{\HOD}{{\rm HOD}}
\newcommand{\Ord}{\mathop{{\rm Ord}}}
\newcommand{\AC}{{\rm AC}}
\newcommand{\DC}{{\rm DC}}
\newcommand{\DCOrd}{\ensuremath{\DC_{{<} \! \Ord}\xspace}}
\newcommand{\CH}{{\rm CH}}
\newcommand{\Coll}{\mathop{\rm Coll}}
\newcommand{\Add}{\mathop{\rm Add}}

\newcommand{\dom}{\text{dom}}
\newcommand{\crit}{\text{crit}}
\newcommand{\la}{\langle}
\newcommand{\ra}{\rangle}
\newcommand{\one}{\mathop{1\hskip-2.5pt {\rm l}}}
\newcommand{\tail}{\text{tail}}
\newcommand{\forces}{\Vdash}
\newcommand{\Ione}{\ensuremath{\rm I_1}\xspace}

\newcommand{\p}{\mathbb P}
\newcommand{\q}{\mathbb Q}
\newcommand{\J}{\mathbb J}
\newcommand{\C}{\mathcal C}

\newcommand{\Los}{\L o\'s}
\newcommand{\Godel}{G\"{o}del}

\newcommand{\unionmodel}[2]{\ensuremath{{#1}_{[#2]}}}
\newcommand{\suppprod}[3]{\ensuremath{{\displaystyle \prod_{\textsubscript{\phantom{A}}{#2}_{\phantom{A}}}}^{(#3)} #1}}

\DeclareMathOperator{\restrict}{\upharpoonright}
\newcommand{\image}{\textnormal{''}}
\newcommand{\treeauto}[1]{\tilde{#1}}
\newcommand{\jextended}{j^+}

\title{ZFC without Power Set II: Reflection Strikes Back}

\author{Victoria Gitman}
\address[V. Gitman]{The City University of New York, CUNY Graduate Center, Mathematics Program, 365 Fifth Avenue, New York, NY 10016}
\email{vgitman@nylogic.org}
\urladdr{https://victoriagitman.github.io/}

\author{Richard Matthews}
\address[R. Matthews] {Univ. Paris Est Cr\'{e}teil, LACL, F-94010}
\email{richard.matthews@u-pec.fr}

\begin{document}

\begin{abstract}
The theory $\ZFC$ implies the scheme that for every cardinal $\delta$ we can make $\delta$ many dependent choices over any definable relation without terminal nodes. Friedman, the first author, and Kanovei constructed a model of $\ZFCminus$ ($\ZFC$ without power set) with largest cardinal $\omega$ in which this principle fails for $\omega$ many choices. In this article we study failures of dependent choice principles over $\ZFCminus$ by considering the notion of big proper classes. A proper class is said to be big if it surjects onto every non-zero ordinal.

We shall see that if one assumes the scheme of dependent choices of any arbitrary set length then every proper class is indeed big. However, by building on work of Zarach, we provide a general framework for separating dependent choice schemes of various lengths by producing models of $\ZFCminus$ with proper classes that are not big. Using a similar idea, we then extend the earlier result by producing a model of $\ZFCminus$ in which there are unboundedly many cardinals but the scheme of dependent choices of length $\omega$ still fails.

Finally, the second author has proven that a model of $\ZFCminus$ cannot have a non-trivial, cofinal, elementary self-embedding for which the von Neumann hierarchy exists up to its critical point. We answer a related question posed by the second author by showing that the existence of such an embedding need not imply the existence of any non-trivial fragment of the von Neumann hierarchy. In particular, that in such a situation $\mathcal{P}(\omega)$ can be a proper class.
\end{abstract}

\maketitle

\section{Introduction}
\noindent Many natural set-theoretic structures satisfy all the axioms of $\ZFC$ excluding the power set axiom. These include the structures $H_{\kappa^+}$ (the collection of all sets whose transitive closure has size at most $\kappa$, where $\kappa$ is a cardinal), forcing extensions of models of $\ZFC$ by pretame (but not tame) class forcing, and first-order structures bi-interpretable with models of the strong second-order set theory Kelley-Morse together with the choice scheme. The set theory that these structures satisfy is the theory $\ZFCminus$, whose axioms consist of the axioms of $\ZFC$ with the collection scheme in place of the replacement scheme and with the well-ordering principle (the assertion that every set can be well-ordered) in place of the axiom of choice (the assertion that every non-empty family of sets has a choice function). The reason for the particular choice of axioms comprising $\ZFCminus$ is that without the existence of power sets we lose certain equivalences between set theoretic assertions that we tend to take for granted.

\begin{definition} \phantom{a}
\begin{itemize}
\item Let $\ZFminusrep$ be the theory $\ZF$ with the power set axiom removed. That is, $\ZFminusrep$ consists of the axioms: extensionality, empty set, pairing, unions, infinity, the foundation scheme, the separation scheme and the replacement scheme.
\item Let $\ZFCminusrep$ denote the theory $\ZFminusrep$ plus the well-ordering principle.
\item Let $\rm{ZF(C)}^{\uminus}$ denote the theory $\rm{ZF(C)}\cminus$ plus the collection scheme.
\end{itemize}
\end{definition}

\noindent Szczepaniak showed that the axiom of choice is not equivalent to the well-ordering principle over $\ZFminus$ (see \cite{zarach82}) and therefore we choose to take the stronger principle when formulating the theory $\ZFCminus$. Zarach showed that the theory $\ZFCminusrep$ does not imply the collection scheme \cite{zarach96}. The first author et al. showed in \cite{gitmanhamkinsjohnstone16} that the theory $\ZFCminusrep$ has many other undesirable behaviors: there are models of $\ZFCminusrep$ in which $\omega_1$ is singular, in which every set of reals is countable but $\omega_1$ exists, and in which the \Los\ Theorem fails for (class) ultrapowers.

Although the theory $\ZFCminus$ avoids these pathological behaviors, there is a number of useful properties of models of full $\ZFC$ that fail or are not known to hold in models of $\ZFCminus$, mostly as a consequence of the absence in these models of a hierarchy akin to the von Neumann hierarchy. It is known that ground model definability, the assertion that the model is definable in its set forcing extensions, can fail in models of $\ZFCminus$ \cite{gitmanjohnstone:groundmodels}. The intermediate model theorem, the assertion that any intermediate model between the model and its set-forcing extension is also its set-forcing extension, can fail \cite{AntosFriedmanGitman:booleanclassforcing}. If there is a non-trivial elementary embedding $j \colon V_{\lambda+1}\to V_{\lambda+1}$, namely the large cardinal axiom \Ione holds, then it gives rise to an elementary embedding $\jextended \colon H_{\lambda^+}\to H_{\lambda^+}$ which witnesses that Kunen's Inconsistency can fail for models of $\ZFCminus$ \cite{matthews20}. It is an open question whether $\HOD$, the collection of all hereditarily ordinal definable sets, is definable in models of $\ZFCminus$.

One of the main themes of this article is the various ways in which the scheme version of dependent choice can fail in models of $\ZFCminus$.

\begin{definition}
The $\DC_\delta$-scheme, for an infinite cardinal $\delta$, asserts for every formula $\varphi(x,y,a)$ that if for every set $x$, there is a set $y$ such that $\varphi(x,y,a)$ holds, then there is a function $f$ on $\delta$ such that for every $\xi<\delta$, $\varphi(f\restrict\xi,f(\xi),a)$ holds.

The $\DCOrd$-scheme is the scheme asserting that the $\DC_\delta$-scheme holds for every cardinal $\delta$.
\end{definition}

\noindent In other words, the $\DC_\delta$ schemes states that we can make $\delta$-many dependent choices along any definable relation without terminal nodes. The $\DC_\delta$-scheme generalizes the dependent choice axiom $\DC_\delta$ which makes the analogous assertion for set relations. The $\DCOrd$-scheme follows from $\ZFC$ by reflecting the definable relation in question to some $V_\alpha$, and then using a well-ordering of $V_\alpha$ to obtain the sequence of dependent choices. It follows that the $\DCOrd$-scheme holds in every structure $H_{\kappa^+}$. It is not known whether pretame class forcing over models of $\ZFC$ preserves the $\DC_{\delta}$-schemes, unless the forcing has no proper class-sized antichains (see \myref{prop:pretamePreservesDC}).

The $\DC_\delta$-schemes have numerous applications. Over $\ZFCminus$, the $\DC_{\omega}$-scheme is equivalent to the \emph{reflection principle} which is the assertion that every formula reflects to a transitive set \cite{friedmangitmankanovei19} (see \myref{th:DC_omegaReflection}). Although, there is no known reformulation of the $\DC_\delta$-scheme for uncountable $\delta$ in terms of a reflection principle, such a reformulation exists under mild existence of power set assumptions (see \myref{DC_deltaReflection}). Over $\ZFCminus$, for a regular cardinal $\delta$, the $\DC_\delta$-scheme implies that every proper class surjects onto $\delta$ (see \myref{prop:bigClasses}). It is not difficult to see that in a model of $\ZFC$ the class partial order $\Add(\Ord,1)$, whose conditions are partial functions from a set of ordinals into $2$ orderd by extension, forces a global well-order without adding sets. This is because, using $\AC$, any set can be coded as a subset of an ordinal and, by genericity, this subset will appear somewhere in the generic class function from $\Ord$ into $2$. We can then well-order the sets by comparing the least location in the generic function where a code appears. In a model of $\ZFCminus$, the forcing $\Add(\Ord, 1)$ is pretame if and only if the $\DCOrd$-scheme holds. Thus, in a model of $\ZFCminus + \DCOrd$-scheme, we can force a global well-order without adding sets, and conversely if we can force a global well-order without adding sets using some forcing, then $\Add(\Ord, 1)$ must be pretame. We will see another application of the $\DCOrd$-scheme shortly to establishing a form of Kunen's Inconsistency for models of $\ZFCminus$.

Friedman et al. \cite{friedmangitmankanovei19} showed that the $\DC_\omega$-scheme can fail in a model of $\ZFCminus$. Moreover, this failure is witnessed by a $\Pi^1_2$ formula, which turns out to be the simplest complexity for which such a failure can occur (see \cite{friedmangitmankanovei19} and Theorem VII.9.2 of \cite{simpsonSOA} for more details). The counterexample model is the $H_{\omega_1}$ of a symmetric submodel of a forcing extension by the iteration of Jensen's forcing along the tree $\omega_1^{{<}\omega}$ (see \myref{subsec:Jensen} for details on Jensen's forcing and this result), in particular, $\omega$ is the largest cardinal in this model. The symmetric submodel in question satisfies $\AC_\omega$, but has a $\Pi^1_2$-definable failure of $\DC$, which translates to its $H_{\omega_1}$ having the requisite properties. 

There are two principle difficulties in constructing such consistency results; having second-order definable failures of $\DC$ and satisfying the full axiom of choice. For example, it is an old result of Jensen that it is possible to produce models of $\ZF$ in which the axiom of choice for families of size at most $\delta$ holds (where $\delta$ is an arbitrary regular cardinal), but $\DC_\omega$ already fails. Furthermore, by Pincus, for any regular cardinal $\delta$ there is a model of $\ZF + \DC_\delta + \neg \DC_{\delta^+}$. We refer the reader to Chapter 8 of Jech's book on the axiom of choice, \cite{jechchoice}, for more details.

In this article we obtain the following failures of the various $\DC_\delta$-schemes in models of $\ZFCminus$.

\begin{duplicate}[\ref{th:DC_omegaHoldsDC_omega_2Fails}]
Suppose that $V\models\ZFC+\CH$. Then every Cohen forcing extension of $V$ has a proper class transitive submodel satisfying $\ZFCminus$ in which the $\DC_{\omega}$-scheme holds, but the $\DC_{\omega_2}$-scheme fails. If we assume further that $V = L$ holds then the $\DC_{\omega_1}$-scheme additionally fails.
\end{duplicate}

\noindent The model above was constructed by Zarach in \cite{zarach82}, and the failure of $\DC_{\omega_1}$ follows by a result of Blass on the Cohen forcing $\Add(\omega,1)$ (see \myref{th:blassCohenReals}). Note that this model, unlike the counterexample model of \cite{friedmangitmankanovei19}, must have unboundedly many cardinals by virtue of being a proper class transitive submodel of a model of $\ZFC$. The second model is constructed by generalizing Zarach's construction as well as generalizing Blass's result to the forcing $\Add(\delta,1)$ (see \myref{th:blassCohenSubsets}).

\begin{duplicate}[\ref{th:DC_deltaHOldsDC_delta++Fails}]
Suppose that $V \models \ZFC + 2^\delta = \delta^+$ for some regular cardinal $\delta$. \linebreak[3] Then every forcing extension of $V$ by the poset $\Add(\delta,1)$ has a proper class transitive submodel satisfying $\ZFCminus$ in which the $\DC_\delta$-scheme holds, but the $\DC_{\delta^{++}}$-scheme fails. If we assume further that $V = L$ holds then the $\DC_{\delta^+}$-scheme also fails.
\end{duplicate}

\noindent Using the idea of union models we extend the result of \cite{friedmangitmankanovei19} to obtain a model of $\ZFCminus$ in which the reflection principle fails and for which there are unboundedly many cardinals.

\begin{duplicate}[\ref{th:WJensenZFC-}, \ref{th:WJensennotReflection}]
Every forcing extension of $L$ by the iteration of Jensen's forcing along the class tree $\Ord^{{<}\omega}$ has a proper class transitive submodel $N$ satisfying $\ZFCminus$, with unboundedly many cardinals, in which the $\DC_{\omega}$-scheme fails.
\end{duplicate}

\noindent The results of this article were originally motivated by a question from the work of the second author on Kunen's Inconsistency in models of $\ZFCminus$ \cite{matthews20}. Suppose that $W\models\ZFCminus$ and $A\subseteq W$. We will say that $W\models\ZFCminus_A$ if $W$ \linebreak[3] continues to satisfy $\ZFCminus$ in the language expanded by a predicate for $A$.

\begin{theorem}[\cite{matthews20}]\label{th:KI}
Suppose that $W\models\ZFCminus$. There is no non-trivial, cofinal, $\Sigma_0$-elementary embedding $j \colon W\to W$ such that $V_{\crit(j)}$ exists in $W$ and $W\models\ZFCminus_j$.
\end{theorem}

\noindent Thus, in particular, the elementary embedding $\jextended \colon H_{\lambda^+}\to H_{\lambda^+}$ resulting from an \Ione-embedding $j \colon V_{\lambda+1}\to V_{\lambda+1}$ cannot be cofinal. Moreover, the second author showed that if the model $W$ additionally satisfies the $\DCOrd$-scheme, then the existence of $V_{\crit(j)}$ follows from the other assumptions \cite{matthews20}. Thus, we have:

\begin{theorem}[\cite{matthews20}]\label{th:criticalPointHierarchy}
Suppose that $W \models \ZFCminus + \DCOrd$-scheme. There is no non-trivial, cofinal, $\Sigma_0$-elementary embedding $j \colon W\to W$ such that $W\models\ZFCminus_j$.
\end{theorem}

\noindent As a natural next step, the second author asked whether the existence of $V_{\crit(j)}$ is truly necessary for \myref{th:KI},

\begin{question}
Is the following situation consistent: There is a non-trivial, cofinal, elementary embedding $j \colon W \rightarrow W$ such that $W \models \ZFCminus_j$?
\end{question}

\noindent In private communications with the second author, Yair Hayut has shown that the above situation is inconsistent, that is there are no non-trivial, cofinal, elementary embeddings $j \colon W \rightarrow W$ for which $W \models \ZFCminus_j$. However, as an initial attempt to answering this question, the second author asked:

\begin{question}
Suppose that $W\models\ZFCminus_j$ for a non-trivial, cofinal $j \colon W\to M$ with $M\subseteq W$. Does $V_{\crit(j)}$ exist in $W$, does $\mathcal \mathcal P (\omega)$ exist in $W$?
\end{question}

\noindent We answer the second question negatively here using models of $\ZFCminus$ in which the $\DC_\delta$-scheme fails for some $\delta$.

\begin{duplicate}[\ref{th:EmbeddingWithoutPowerSet}]
There is a model $W\models\ZFCminus$ in which $\mathcal P (\omega)$ does not exist and which has a definable, non-trivial, cofinal, elementary embedding $j \colon W\to M\subseteq W$.
\end{duplicate}

\begin{remark}[A note on the title]
This paper should be seen as a continuation of the study of models of set theory without power set carried out in \cite{gitmanhamkinsjohnstone16}. In particular, we see some of the pathological properties that can arise in such models when we don't assume the $\DCOrd$-scheme. As such, we have titled this work as part two on the question of what is $\ZFC$ without power set.
\end{remark}

\section{Preliminaries}\label{sec:prelim}
\subsection{Reflection}

Although, the connection between the $\DC_\delta$-schemes and reflection is not necessary for any of the arguments in this article, we nevertheless devote this section to exploring this unexpected connection. The connection also lends to the title of the paper.

\begin{theorem}\label{th:DC_omegaReflection}
Over $\ZFCminus$, the $\DC_\omega$-scheme is equivalent to the scheme asserting for every formula $\psi(\vec x,a)$ that there is a transitive set $M$ with $a\in M$ reflecting $\psi(\vec x,a)$.
\end{theorem}

\noindent We sketch the proof of this result, which has appeared in \cite{friedmangitmankanovei19}.

\begin{proof}
First, observe that the $\DC_{\omega}$-scheme is equivalent to the $\DC^*_\omega$-scheme asserting that for every definable relation $\varphi(x,y,a)$ without terminal nodes, there is a sequence $\la b_i\mid i<\omega\ra$ such that $\varphi(b_i,b_{i+1},a)$ holds for every $i<\omega$.

Suppose that $W$ satisfies $\ZFCminus$ and the reflection assertion. Suppose that $\varphi(x,y,a)$ is a relation without terminal nodes. Let $M$ be a transitive set with $a\in M$ which reflects $\varphi(x,y,a)$ together with the assertion that $\varphi(x,y,a)$ has no terminal nodes. Since $M$ is a set, there is a well-ordering $w$ of $M$ in $W$. Let $b_0$ be the least element of $M$ according to $w$. Let $b_1$ be the least element of $M$ according to $w$ such that $M\models\varphi(b_0,b_1, a)$, which exists since $M$ knows that $\varphi(x,y,a)$ has no terminal nodes. Given that we have chosen $b_n\in M$, let $b_{n+1}$ be the least element of $M$ according to $w$ such that $M\models\varphi(b_n,b_{n+1},a)$. Clearly, the sequence $\la b_i\mid i<\omega\ra$ witnesses the $\DC_{\omega}^*$-scheme for $\varphi(x,y,a)$.

Next, suppose that $W \models \ZFCminus + \DC_\omega$-scheme. The result will follow from induction on formulas with the only critical case being those formulas of the form $\exists x \psi(x, u)$. So suppose that the statement has been proven for $\psi(x, u)$. Observe that, by collection and the induction hypothesis, for any set $A$ there is a transitive set $A^\psi$ containing $A$ which reflects $\psi(x, u)$ and such that
\[
\forall u \in A \, \exists x \, \psi(x, u) \longrightarrow \forall u \in A \, \exists x \in A^\psi \, \psi(x, u).
\]
Fix a set $a$. Let the formula $\varphi(x, y, a)$ assert that whenever $x$ is a sequence of some finite length $n$ such that $x_0 = \{ a \}$, and $x_{i+1} = x_i^\psi$, then $y = x_{n-1}^\psi$. By the above argument, this relation has no terminal nodes. Using the $\DC_\omega$-scheme, the union of an $\omega$-sequence of dependent choices along $\varphi(x,y,a)$ is a transitive set reflecting $\exists x \psi(x,u)$ and containing $a$.
\end{proof}

\noindent It is worth noting why the above argument does not also show that the reflection principle implies the $\DC_\delta$-scheme for uncountable cardinals. The issue is that the transitive set $M$ need not be closed under infinite sequences which are elements of $W$. Therefore, if we reflect our formula to some arbitrary set $M$ our attempt to externally choose the $b_\alpha$ may fail because we cannot ensure at limit stages that our collection of previous choices forms a set in $M$.

Unfortunately, we do not know whether there is reformulation of the $\DC_{\delta}$-scheme for uncountable $\delta$ in terms of some reflecting principle. However, we do have the following weaker result.

\begin{theorem}\label{DC_deltaReflection}
Suppose that $W$ satisfies $\ZFCminus$ and $\delta$ is a regular cardinal in $W$ such that $\gamma^{{<}\delta}$ exists for every cardinal $\gamma$. Then in $W$, the $\DC_\delta$-scheme holds if and only if for every formula $\psi(\vec x,a)$, there is a transitive set $M$ with $a\in M$ and $M^{{<}\delta}\subseteq M$ reflecting $\psi(\vec x,a)$.
\end{theorem}

\begin{proof}
Suppose that $W$ satisfies the reflection assertion. Fix a definable relation $\varphi(x,y,a)$ without terminal nodes and let $M$ be a transitive set with $a\in M$ and $M^{{<}\delta}\subseteq M$ which reflects $\varphi(x,y,a)$ and the assertion that this relation has no terminal nodes. We construct a sequence of $\delta$-many dependent choices as in the proof of \myref{th:DC_omegaReflection}, using the closure of $M$ to get through the limit stages in the construction.

Suppose next that the $\DC_\delta$-scheme holds in $W$.
We will say that a set $A^\psi$ is a \emph{$\delta$-transitive closure of a set $A$ for a formula $\psi$} if it is a transitive set containing $A$ which reflects $\psi$, is closed under existential witnesses for $\psi$ from $W$, and closed under ${<}\delta$-sequences. We need the assumption that $\gamma^{{<}\delta}$ exists for every cardinal $\gamma$ to ensure that every set can be closed under ${<}\delta$-sequences. From here the argument proceeds exactly as in the proof of \myref{th:DC_omegaReflection}.
\end{proof}

\subsection{Big classes}

Given a cardinal $\delta$, let us say that a class is \emph{$\delta$-big} if it surjects onto $\delta$. We will say that a class is simply \emph{big} if it surjects onto every cardinal. It is easy to see that proper classes don't need to be big in weak set theories. For example, consider the model $L_{\aleph_\omega^L}$, which satisfies Kripke-Platek set theory, ${\rm KP}$. The cardinals of $L_{\aleph_\omega^L}$ is a proper class from the point of view of this model, but this class obviously cannot surject onto $\aleph_1^L$ because externally we know that it is countable. We will see in Sections \myref{sec:Zarach} and \myref{sec:GeneralizedZarach} that proper classes do not need to be big in models of $\ZFCminus$ either. However, $\ZFCminus+\DCOrd$-scheme implies that every proper class is big.
\begin{proposition}\label{prop:bigClasses}
In a model of $\ZFCminus$, the $\DC_\delta$-scheme implies that every proper class is $\delta$-big. It follows that over $\ZFCminus$, the $\DCOrd$-scheme implies that every proper class is big.
\end{proposition}
\begin{proof}
Let $W\models\ZFCminus+\DC_{\delta}$-scheme for some regular cardinal $\delta$. Consider a proper class $\mathcal A$ defined by a formula $\psi(x,a)$. Let $\varphi(x,y,a)$ be a formula asserting that whenever $x$ is a function on an ordinal $\xi$ such that $x(\eta)\in \mathcal A$ for all $\eta<\xi$, then $y\in \mathcal A$ and $y\neq x(\eta)$ for any $\eta<\xi$. Since $\mathcal A$ is a proper class, the relation $\varphi(x,y,a)$ has no terminal nodes. Thus, by the $\DC_\delta$-scheme, there is a function $f$ on $\delta$ such that for all $\xi<\delta$, $\varphi(f\restrict \xi,f(\xi),a)$ holds. The function $f$ gives a subset of $\mathcal A$ of cardinality $\delta$.
\end{proof}

\noindent The model constructed in \cite{friedmangitmankanovei19} to show that the $\DC_\omega$-scheme can fail also shows that the converse to \myref{prop:bigClasses} does not hold. Recall that the model is the $H_{\omega_1}$ of a model $W\models\ZF+\AC_{\omega}$. But $\AC_\omega$ implies that every set surjects onto $\omega$, the largest cardinal of the model. Thus, in that model, every class is big. In \myref{sec:LargeModelNegDC}, we will strengthen this by showing that the $\DC_\omega$-scheme can fail in a model of $\ZFCminus$ with unboundedly many cardinals in which every proper class is big. On the other hand, adding small proper classes will be one of our main tools in this article for constructing models of $\ZFCminus$ with various violations of the $\DC_\delta$-scheme.

\myref{th:criticalPointHierarchy} from the introduction was obtained by showing that whenever we have a model $W\models\ZFCminus$, in which every proper class is big and $j \colon W\to M\subseteq W$ is an elementary embedding with a critical point, $V_{\crit(j)}$ exists in $W$ \cite{matthews20}. We will quickly reprove the theorem here to emphasize the exact assumptions and demonstrate how big classes are used in the proof.

\begin{theorem}\label{th:elementaryEmbeddingHierarchy}
Suppose that $W\models\ZFCminus$ and every proper class is big in $W$. If $j \colon W\to M\subseteq W$ is an elementary embedding with critical point $\kappa$, then $V_\kappa$ exists in $W$.
\end{theorem}

\noindent Note that we are not assuming that $W\models\ZFCminus_j$ or that $j$ is cofinal.

\begin{proof}
First, observe that if $\alpha<\kappa$ and $A\subseteq \alpha$, then $j(A)=A$. Next, let's argue that $\mathcal P (\alpha)$ exists for every $\alpha<\kappa$. Fix $\alpha<\kappa$. Suppose towards a contradiction that $\mathcal P (\alpha)$ is a proper class. Then by our assumption that every proper class is big, there is a surjection from $\mathcal P (\alpha)$ onto $\kappa$. Applying collection, we can obtain a set $B\subseteq \mathcal P (\alpha)$ for which there is a surjection $h \colon B\to \kappa$. By elementarity, $j(h) \colon j(B)\to j(\kappa)$ is a surjection onto $j(\kappa)$. Observe that $b\in B$ if and only if $b=j(b)\in j(B)$, and hence $B=j(B)$. Thus, $j(h) \colon B \to j(\kappa)$. Also, for every $b\in B$, $$j(h)(b)=j(h)(j(b))=j(h(b))=h(b).$$ Thus, the range of $j(h)$ is $\kappa$ contradicting that $j(h)$ is a surjection onto $j(\kappa)$. Now, a standard argument shows that $|\mathcal P (\alpha)|<\kappa$ for every $\alpha<\kappa$, and that $\kappa$ is regular.

Next, let's argue that $V_\alpha$ exists for every $\alpha\leq\kappa$. Suppose inductively that we have shown that $V_\alpha$ exists and $|V_\alpha|=\beta<\kappa$. Then we can use a bijection $f \colon \beta\to V_\alpha$ and the previously shown fact that $\mathcal P (\beta)$ exists and $|\mathcal P (\beta)|<\kappa$ to argue that $V_{\alpha+1}$ exists and $|V_{\alpha+1}|<\kappa$. We then use collection to argue that $V_\lambda$ exist for limit $\lambda$, and use the regularity of $\kappa$ to argue that $|V_\lambda|<\kappa$ for $\lambda<\kappa$.
\end{proof}

\subsection{Class forcing}

In \myref{sec:LargeModelNegDC}, we use class forcing to construct a model of $\ZFCminus$ with unboundedly many cardinals and all big proper classes in which the $\DC_\omega$-scheme fails. Here, we briefly summarize the relevant properties of class forcing which we shall use in that argument.

Class forcing is best interpreted when working over a model of some second-order set theory. Second-order set theory is formalized in a two-sorted logic with separate sorts (variables and quantifiers) for sets and classes. Thus, unlike in first-order set theory, in this setting classes are actual elements of the model and not just objects of the meta-theory. Models of second-order set theory are triples $\mathcal W =\la W,\in,\mathcal C\ra$ where $W$ is the sets of the model, $\mathcal C$ is the classes, and $\in$ is the membership relation between sets, as well as between sets and classes, letting us know of which sets each of the classes is composed. Let $\GBminus$ denote the second-order set theory whose axioms for sets are $\ZFminus$ and whose axioms for classes consist of extensionality, the class collection axiom asserting that for every class relation whose domain is restricted to a set, there is a set of witnesses of the relation's image, and the first-order comprehension scheme asserting that every first-order definable collection of sets is a class. Furthermore, we let $\GBcminus$ be $\GBminus$ plus the axiom of choice and $\GBCminus$ be $\GBminus$ with the global well-order axiom, which asserts that there is a bijection between $W$ and $\Ord$. By replacing $\ZFminus$ with $\ZF$ or $\ZFC$ in the theory $\GBminus$, we obtain the \Godel-Bernays set theories of $\GB$ and $\GBc$ respectively. Every model of $\ZFCminus$ with a definable global well-order is naturally a model of $\GBCminus$ and every model of $\ZFC$ with a definable global well-order is naturally a model of $\GBC$. By forcing with $\Add(\Ord,1)$, we can show that every model of $\ZFC$ has a class forcing extension with the same sets and a global well-order. Thus, every model of $\ZFC$ has a class forcing extension with the same sets that is a model of $\GBC$. As explained in the introduction, the analogous fact is true only for models $\ZFCminus$ provided that the $\DCOrd$-scheme holds.

A class forcing notion in a model $\mathcal W=\la W,\in,\mathcal C\ra\models\GBminus$ is a partial order $\p\in \mathcal C$. Suppose that $\mathcal W=\la W,\in,\mathcal C\ra\models\GBminus$ and $\p\in\mathcal C$ is class forcing notion. We say that $G\subseteq \p$ is $\mathcal W$-\emph{generic} if $G$ meets every dense class $D\subseteq \p$ from $\mathcal C$. The forcing extension is $\mathcal W[G]=\la W[G],\in,\mathcal C[G]\ra$, where $W[G]$ is the collection of all interpretations of (the usual) $\p$-names by $G$, and $\mathcal C[G]$ is the collection of all interpretations of the class $\p$-names by $G$, where a \emph{class $\p$-name} is a class whose elements are pairs $\la \dot x ,p\ra$ where $\dot x$ is a $\p$-name and $p\in\p$.

In a number of significant ways, class forcing does not behave as nicely as set forcing. It is easy to see, for example by forcing with $\Coll(\omega,\Ord)$ (conditions are finite functions from $\omega$ to the ordinals ordered by extension) to collapse $\Ord$ to $\omega$, that class forcing need not preserve replacement to the forcing extension. The forcing relations for a class forcing notion need not be definable (or more generally need not be a class). For example, a model of $\GBC$ whose classes are definable collections, can have a class forcing notion for which the forcing relation on atomic formulas is not definable \cite{holykrapflucke16}. However, there is a class of well-behaved class forcing notions, the pretame forcings, which avoid these pathological behaviors.

\begin{definition}
Suppose that $\mathcal W=\la W,\in,\mathcal C\ra\models\GBminus$. A notion of class forcing $\mathbb{P}\in\mathcal C$ is \emph{pretame} if for every $p \in \mathbb{P}$ and any sequence of classes \mbox{$\langle D_i \mid i \in I \rangle \in \mathcal{C}$}, with $I \in W$, such that each $D_i$ is dense below $p$, there is a condition $q \leq p$ and a sequence $\langle d_i \mid i \in I \rangle \in W$ such that for every $i \in I$, $d_i \subseteq D_i$ and $d_i$ is predense below $q$.
\end{definition}

\begin{theorem}\phantom{a}
\begin{enumerate}
\item $($Friedman \cite{friedmanclassforcing}$)$ Pretame class forcing notions preserve $\GBminus$ to the forcing extension.
\item $($Stanley \cite{holykrapfschlicht18}$)$ Pretame class forcing notions have definable forcing relations.
\end{enumerate}
\end{theorem}

\noindent In the context of models of second-order set theory, let's redefine the $\DC_\delta$-scheme to assert that we can make $\delta$-many dependent choices over every class (not just definable) relation without terminal nodes. In particular, all our results will follow for models in which the only classes are the definable collections. Although, we will not make use of the following proposition in the rest of the article, the result fits into our analysis of the $\DCOrd$-scheme. However, this result will require that our class forcing satisfies an additional assumption known as the \emph{Maximality Principle}. By \cite{holykrapfschlicht18}, over $\GBCminus$, this is known to be equivalent to the assumption that every anti-chain is a set.

\begin{definition}
Suppose that $\mathcal W = \la W, \in, \mathcal C \ra \models \GBminus$. A notion of class forcing $\mathbb P \in \mathcal C$ satisfies the \emph{Maximality Principle} if whenever $p \Vdash \exists x \varphi(x, \dot y, \dot \Gamma)$ for some $p \in \mathbb P$ and formula $\varphi$ with class name parameter $\dot \Gamma \in \mathcal C$ and set name parameter $\dot y \in W^{\mathbb P}$, then there exists some $\dot a \in W^{\mathbb P}$ such that $p \Vdash \varphi(\dot a, \dot y, \dot \Gamma)$.
\end{definition}

\begin{proposition}\label{prop:pretamePreservesDC}
Suppose that $\mathcal W=\la W,\in,\mathcal C\ra\models\GBminus+\DC_\delta$-scheme for some regular cardinal $\delta$. Then every pretame forcing extension of $\mathcal W$ which satisfies the Maximality Principle and in which $\delta$ remains regular satisfies the $\DC_\delta$-scheme.
\end{proposition}

\begin{proof}
Suppose that $\p\in\mathcal C$ is a pretame forcing notion. Let $G\subseteq \p$ be $\mathcal W$-generic. Let $R\in\C[G]$ be a class relation without terminal nodes. Let $\dot R$ be a class $\p$-name for $R$ and let $p\in\p$ be a condition forcing that $\dot R$ does not have terminal nodes. Given a sequence $x$ of $\p$-names of length some ordinal $\xi$, let $\dot x^{(\xi)}$ denote the canonical $\p$-name for a sequence of length $\xi$ whose $\eta$-th element, for $\eta<\xi$, is the interpretation of $x(\eta)$. Using the Maximality Principle, let $\varphi(x,\dot y,p,\dot R,\p)$ be a formula asserting, over $\mathcal W$, that whenever $x$ is a sequence of $\p$-names of some ordinal length $\xi$ and $p$ forces that $\dot x^{(\eta)}\,\dot R\, x(\eta)$, then $\dot y$ is a $\p$-name and $p$ forces that $\dot x^{(\xi)}\,\dot R\, \dot y$. Since $p$ forces that $\dot R$ has no terminal nodes, the relation given by $\varphi$ has no terminal nodes either. Thus, by the $\DC_\delta$-scheme, we can make $\delta$-many choices along the relation given by $\varphi$. Let $f$ be the function with domain $\delta$ witnessing this. Then $\dot f^{(\delta)}_G$ witnesses that we can make $\delta$-many dependent choices over the relation $R$.
\end{proof}

\noindent The next proposition gives a useful criterion for pretameness.

\begin{proposition}\label{prop:chainConditionPretame}
Suppose that $\mathcal W=\la W,\in,\mathcal C\ra\models\GBminus+\DC_{\delta}$-scheme, for a regular cardinal $\delta$, and $\p\in\C$ is a class forcing notion with the $\delta$-cc. Then $\p$ is pretame.
\end{proposition}

\begin{proof}
We will argue that every dense class $D\subseteq\p$ has a set maximal antichain contained in it, these antichains will then witness pretameness. Suppose towards a contradiction that there is no set maximal antichain contained in $D$. Let $\varphi(x,y,\p,D)$ be a formula asserting, over $\mathcal W$, that whenever $x$ is a sequence of incompatible elements of $D$ of some ordinal length $\xi$, then $y\in D$ and $y$ is incompatible with all elements of the sequence $x$. Since there is no set maximal antichain contained in $D$, the relation given by $\varphi$ has no terminal nodes. Thus, we can make $\delta$-many dependent choices along it, which contradicts our assumption that $\p$ has the $\delta$-cc.
\end{proof}

\subsection{Jensen's forcing}\label{subsec:Jensen}

In any universe $V\models\ZFC+\diamondsuit$ we can construct a subposet $\J$ of Sacks forcing (elements are perfect trees ordered by the subtree relation) with the following two key properties.

\begin{theorem}[Jensen \cite{jensen70}]\phantom{a}
\begin{enumerate}
\item The poset $\J$ has the ccc.
\item Suppose that $r$ is a $V$-generic real for $\J$. Then in $V[r]$, $r$ is the unique $V$-generic real for $\J$.
\end{enumerate}
\end{theorem}

\noindent Such a poset $\J$ was first constructed by Jensen in $L$ \cite{jensen70}. The choice of the $\diamondsuit$-sequence can potentially yield different such posets $\J$. In $L$, Jensen used the canonical $\diamondsuit$-sequence (defined by taking the least counterexample at each stage) to construct such a poset $\J$ with the additional property that the unique $L[G]$-generic real added by $\J$ is a $\Pi^1_2$-definable singleton \cite{jensen70}. This is the lowest possible such complexity because $\Pi^1_2$-definable singleton reals must be constructible by Shoenfield's absoluteness.

Before we proceed, let us introduce a general notation for the product of $\mu$ many copies of a forcing with support of size less than $\delta$, which we will use throughout this article.

\begin{definition}
Let $\p$ be a forcing notion and let $\delta \leq \mu$ be regular cardinals. Let $\suppprod{\p}{\mu}{\delta}$ denote the product forcing of $\mu$ many copies of $\p$ with ${<}\delta$ support.
\end{definition}

\medskip

\noindent Observe that we can treat conditions in $\suppprod{\p}{\mu}{\delta}$ as functions $f \colon \mu \rightarrow \p$ such that $\{ \xi \in \mu \mid f(\xi) \neq \mathbbm{1} \} < \delta$. Lyubetsky and Kanovei showed that the poset $\suppprod{\J}{\omega}{\omega}$ has the ccc and the following uniqueness of generics property.

\medskip

\vbox{
\begin{theorem}[Lyubetsky, Kanovei \cite{kanovei:productOfJensenReals}]\label{th:JensenUniqueness} Suppose that $V\models\ZFC+\diamondsuit$.
\begin{enumerate}
\item The poset $\suppprod{\J}{\omega}{\omega}$ has the ccc.
\item Suppose that $G\subseteq \suppprod{\J}{\omega}{\omega}$ is $V$-generic. Then in $V[G]$, the $V$-generic reals for $\J$ are precisely the $\omega$-many reals coming from the slices of $G$.
\end{enumerate}
\end{theorem}
}

\noindent In fact, an application of the $\Delta$-system lemma shows that any length finite-support product of $\J$ has these two properties.

We will say that a forcing iteration $\p_n$, of length $n$, is an \emph{iteration of subposets of Sacks forcing} if every initial segment of $\p_n$ forces that the next poset in the iteration is a subposet of the Sacks forcing of that extension. In universes where $\J$ can be constructed, we can construct iterations $\J_n$, for any $n<\omega$, of subposets of Sacks forcing with the following key properties \cite{friedmangitmankanovei19}.
\begin{enumerate}
\item If $m<n$, then $\J_n\restrict m=\J_m$.
\item $\J_n$ has the ccc.
\item Suppose that $\la r_1,\ldots,r_n\ra$ is a $V$-generic sequence of reals for $\J_n$. Then in $V[\la r_1,\ldots,r_n\ra]$, this is the unique $V$-generic sequence of reals for $\J_n$.
\end{enumerate}

Let $\vec \J=\la \J_n\mid n<\omega\ra$. Let $X$ be any set or class and consider the tree $X^{{<\omega}}$ of finite sequences from $X$ ordered by extension. Let $\mathcal T\subseteq X^{{<}\omega}$ be a sub-tree. Let $\p(\vec {\J},\mathcal T)$ be the (possibly class) poset whose elements are functions $f_T$ on a finite subtree $T$ of $\mathcal T$ such that for nodes $s$ on level $n$ of $T$, $f_T(s)\in \J_n$ and for nodes $s<t$ in $T$, we have that $f_T(t)\restrict\text{len}(s)= f_T(s)$. The ordering is given by $f_T\leq g_S$ provided that $T$ extends $S$ and for every node $s\in S$, we have that $f_T(s)\leq g_S(s)$. We call the poset $\p(\vec {\J},\mathcal T)$ , an iteration of Jensen's forcing along the tree $\mathcal T$. It is proven in \cite{friedmangitmankanovei19} that the analogue of the properties in Theorem 2.12 also hold for the tree version of the forcing, which we state in $\GB$ to handle the possibility that $X$ is a class, which would imply that the resulting poset is a class forcing.

\begin{theorem}[\cite{friedmangitmankanovei19}]\label{th:JensenTreeForcing} Suppose that $\mathcal W=\la W,\in,\mathcal C\ra \models \GBc + \diamondsuit$ and $X$ is any set or class.
\begin{enumerate}
\item The poset $\p(\vec {\J},X^{{<}\omega})$ has the ccc.
\item Suppose that $G\subseteq \p(\vec {\J},X^{{<}\omega})$ is $\mathcal W$-generic. Then the $\mathcal W$-generic sequences $\la r_1,\ldots,r_n\ra$ for $\J_n$ in $\mathcal W[G]$ are precisely the sequences added by nodes of $X^{{<}\omega}$ on level $n$.
\end{enumerate}
\end{theorem}

\begin{proposition}\label{prop:JensenTreeForcingPretame}
Suppose that $\mathcal W=\la W,\in,\mathcal C\ra \models \GBc + \diamondsuit$ and $X$ is any set or class. Then the poset $\p(\vec {\J},X^{{<}\omega})$ is pretame.
\end{proposition}

\begin{proof}
This follows by combining \myref{th:JensenTreeForcing} with \myref{prop:chainConditionPretame}.
\end{proof}

\section{Zarach's union models of \texorpdfstring{$\ZFC^-$}{ZFC-}}\label{sec:Zarach}
\noindent In \cite{zarach82}, Zarach gave a general construction for producing interesting models of $\ZFCminus$ as unions of models of $\ZFC$ arising as transitive submodels of a carefully chosen forcing extension. Because of the style in which it was presented, we have decided to rewrite his construction using modern notation.

Suppose that $V\models\ZFC$. Let $\p\in V$ be a poset such that $\q=\suppprod{\p}{\omega}{\omega}$ is isomorphic to $\p$. Let us call an automorphism $\pi$ of $\q$ \emph{coordinate-switching} if there is an automorphism $\treeauto{\pi}$ of $\omega$ such that, for any condition $p$, $\pi(p)=q$, where $q$ is defined by $q(i)=p(\treeauto{\pi}^{-1}(i))$, namely $\pi$ simply switches coordinates according to $\treeauto{\pi}$.

Let $G\subseteq \q$ be $V$-generic. For $n<\omega$, let
\begin{enumerate}
\item $G_n$ be the restriction of $G$ to the first $n$ coordinates,
\item $G_{\{n\}}$ be the restriction of $G$ to the $n$-th coordinate,
\item $G_{n,\tail}$ be the tail of $G$ after $n$ ($G=G_n\times G_{n,\tail}$).
\end{enumerate}
Let $\unionmodel{V}{n} = V[G_n]$. Let $W^V_G=\bigcup_{n<\omega}\unionmodel{V}{n}$.

By uniform ground model definability (\cite{laver:groundmodel}, \cite{woodin:groundmodel}), $W^V_G$ is definable in $V[G]$ from the generic filter $G$, say by the formula $\varphi(x,\p,G)$. Thus, to every formula $\psi(x)$, there corresponds a formula $\psi^W(x,y,z)$ such that for every $a\in W^V_G$, $W^V_G\models\psi(a)$ if and only if $V[G]\models\psi^W(a,\p,G)$. Observe that if $\pi$ is a coordinate-switching automorphism of $\q$, then $\varphi(x,\p,G)$ and $\varphi(x,\p,\pi\image G)$ both define $W^V_G$ in $V[G]$. Hence also, for every formula $\psi(x)$, we have that for every $a\in W^V_G$, $V[G]$ satisfies $\psi^W(a,\p,G)$ if and only if $V[G]$ satisfies $\psi^W(a,\p,\pi\image G)$. Let $\dot G$ be the canonical $\q$-name for the generic filter.

\begin{proposition}\label{prop:OneDecidesSomeWstatements}
If $p\forces\psi^W(\check a,\check{\p},\dot G)$, then $\one\forces\psi^W(\check a,\check{\p},\dot G)$.
\end{proposition}

\begin{proof}
Suppose for a contradiction that $p\forces\psi^W(\check a,\check{\p},\dot G)$ and $q\forces\neg \psi^W(\check a,\check{\p},\dot G)$. Let $n$ be above the domains of $p$ and $q$. Let $\pi$ be a coordinate-switching automorphism that switches the coordinates in the domain of $p$ to some coordinates above $n$. Then $\pi(p)$ and $q$ are clearly compatible, and $\pi(p)\forces\psi^W(\check a,\check{\p},\pi(\dot G))$. Since $\pi(\dot G)_{\pi\image G}=G$ (the image of $\pi\image G$ under the coordinate-switching automorphism $\pi^{-1}$), by our argument above, $\pi(p)\forces\psi^W(\check a,\check{\p},\dot G)$ as well, which is the desired contradiction.
\end{proof}

\noindent In order to show that $W^V_G \models \ZFCminus$, we shall prove that $W^V_G = \bigcup_{n < \omega} \unionmodel{W}{n}$ where $\unionmodel{W}{n} \prec W^V_G$ for each $n < \omega$. To do this, fix an isomorphism $h \colon \p\cong \q=\suppprod{\p}{\omega}{\omega}$. Since each $G_{\{n\}}$ is $V$-generic for $\p$, it follows that $G^{(n)} := h \image G_{\{n\}}$ is $V$-generic for $\q$. As before, given $G^{(n)}$, for each $m < \omega$ we can obtain generics $G^{(n)}_{m}$, $G^{(n)}_{\{m\}}$, and $G^{(n)}_{m,\tail}$. Now, for $n, m < \omega$, let
\[
\unionmodel{(\unionmodel{V}{n})}{m} = V[G_n][G^{(n)}_m] = \unionmodel{V}{n}[G^{(n)}_m]
\]
and let $W^{\unionmodel{V}{n}}_{G^{(n)}}=\bigcup_{m<\omega} \unionmodel{(\unionmodel{V}{n})}{m}$. Thus, we have:
\[
\unionmodel{V}{n} \subseteq \unionmodel{(\unionmodel{V}{n})}{1} \subseteq \unionmodel{(\unionmodel{V}{n})}{2} \subseteq \cdots \subseteq \unionmodel{(\unionmodel{V}{n})}{m} \subseteq \cdots \subseteq \bigcup_{m<\omega}\unionmodel{(\unionmodel{V}{n})}{m} = W^{\unionmodel{V}{n}}_{G^{(n)}} \subseteq \unionmodel{V}{n+1}.
\]
For every $m,n<\omega$, let $H^{(n,m)}$ be the $\unionmodel{(\unionmodel{V}{n})}{m}$-generic filter for $\q$ obtained from $G^{(n)}_{m,\tail}$ via the isomorphism of $\q$ with its tail after $m$. Thus, for every $m,n<\omega$, we have that
\[
\unionmodel{V}{n+1} = \unionmodel{V}{n}[G^{(n)}_m \times G^{(n)}_{m, \tail}] =\unionmodel{(\unionmodel{V}{n})}{m}[H^{(n,m)}]
\]
is a $\q$-forcing extension of $\unionmodel{(\unionmodel{V}{n})}{m}$. Moreover, this yields $W^{\unionmodel{V}{n}}_{G^{(n)}}=W^{\unionmodel{(\unionmodel{V}{n})}{m}}_{H^{(n,m)}}$ in $\unionmodel{V}{n+1}$.

We shall now argue that $W^{\unionmodel{V}{n}}_{G^{(n)}}\prec W^V_G$ for every $n<\omega$. This extremely powerful key lemma will yield many of our desired results.

\begin{lemma}\label{lem:elementarityW}
For every $n<\omega$, $W^{\unionmodel{V}{n}}_{G^{(n)}}\prec W^V_G$.
\end{lemma}

\begin{proof}
Fix a formula $\psi(x)$ and $a\in W^{\unionmodel{V}{n}}_{G^{(n)}}$ such that $W^{\unionmodel{V}{n}}_{G^{(n)}}\models\psi(a)$. Since $W^{\unionmodel{V}{n}}_{G^{(n)}}$ is the union of $\unionmodel{(\unionmodel{V}{n})}{m}$ for $m < \omega$, let $m$ be such that $a\in \unionmodel{(\unionmodel{V}{n})}{m}$. Thus, by \myref{prop:OneDecidesSomeWstatements}, $\one\forces \psi^W(\check a,\check{\p},\dot G)$ over $\unionmodel{(\unionmodel{V}{n})}{m}$. By the above argument, we have $\unionmodel{(\unionmodel{V}{n})}{m}[H^{(n,m)}]=\unionmodel{V}{n+1}$, and thus
\[
\unionmodel{(\unionmodel{V}{n})}{m}[H^{(n,m)}][G_{n,\tail}]=V[G].
\]
Via the obvious isomorphism of $\p\times\q$ and $\q$, we can view $H^{(n,m)}*G_{n,\tail}$ as a $\unionmodel{(\unionmodel{V}{n})}{m}$-generic for $\q$ with $H^{(n,m)}$ being the generic on the first coordinate. Consider the model $W^{\unionmodel{(\unionmodel{V}{n})}{m}}_{H^{(n,m)}*G_{n,\tail}}$ obtained from the generic $H^{(n,m)}*G_{n,\tail}$. We have

\begin{eqnarray*}
W^{\unionmodel{(\unionmodel{V}{n})}{m}}_{H^{(n,m)}*G_{n,\tail}} & = & \unionmodel{(\unionmodel{V}{n})}{m}[H^{(n,m)}] \cup \bigcup_{i < \omega} \unionmodel{(\unionmodel{V}{n})}{m}[H^{(n, m)}][G_{n + 1 + i}] \\
& = & \bigcup_{i < \omega} \unionmodel{V}{n + i} \\
& = & W^V_G.
\end{eqnarray*}

\noindent But we already showed that $\unionmodel{(\unionmodel{V}{n})}{m}$ satisfies that $\one\forces \psi^W(\check a,\check{\p},\dot G)$. Thus, \linebreak[4] \mbox{$W^V_G\models\psi(a)$}.
\end{proof}

\noindent Next, let us see what theory the model $W^V_G$ satisfies. Observe right away that $W^V_G$ cannot be a model of the power set axiom because $\mathcal P (\p)$ does not exist in $W^V_G$.

\begin{theorem}[Zarach \cite{zarach96}]
$W^V_G\models\ZFCminus$.
\end{theorem}

\begin{proof}
It is clear that $W^V_G$ satisfies extensionality, empty set, pairing, unions, infinity, and the foundation scheme. Also, $W^V_G$ satisfies the well-ordering principle because any set in $W^V_G$ is in some $\unionmodel{V}{n}\models\ZFC$. So it remains to argue that $W^V_G$ satisfies the separation and collection schemes. First, let's do separation. Fix a formula $\psi(x,y)$ and some $a,b\in W^V_G$. We need to argue that
\[
c= \big\{ x\in a\mid W^V_G\models \psi(x,b) \big\}
\]
is in $W^V_G$. Let $n$ be large enough so that $a,b\in \unionmodel{V}{n}$. By \myref{lem:elementarityW}, $W^{\unionmodel{V}{n}}_{G^{(n)}}\prec W^V_G$. Thus,
\[
c= \big\{ x\in a\mid W^{\unionmodel{V}{n}}_{G^{(n)}}\models \psi(x,b) \big\}.
\]
But $W^{\unionmodel{V}{n}}_{G^{(n)}}$ is a definable submodel of $\unionmodel{V}{n+1}\subseteq W^V_G$, and $\unionmodel{V}{n+1}\models\ZFC$. Thus, $c\in \unionmodel{V}{n+1}$ by separation, and hence $c\in W^V_G$.

Next, let's do collection. Suppose that $\psi(x,y,z)$ is a formula and $a,b\in W^V_G$ such that
\[
W^V_G\models\forall x\in a\,\exists y\,\psi(x,y,b).
\]
Let $n$ be large enough so that $a,b\in \unionmodel{V}{n}$. By \myref{lem:elementarityW}, $W^{\unionmodel{V}{n}}_{G^{(n)}}\prec W^V_G$, and so $W^{\unionmodel{V}{n}}_{G^{(n)}}\models\forall x\in a\,\exists y\,\psi(x,y,b)$. It follows that $\unionmodel{V}{n+1}$ has a collecting set for $\psi(x,y,b)$, and hence so does $W^V_G$.
\end{proof}

\begin{theorem}[Zarach \cite{zarach96}]
$W^V_G\models\DC_\omega$-scheme.
\end{theorem}

\begin{proof}
Fix a formula $\psi(x,y,a)$ defining over $W^V_G$ a relation without terminal nodes. Let $n$ be large enough so that $a\in \unionmodel{V}{n}$. By \myref{lem:elementarityW}, $W^{\unionmodel{V}{n}}_{G^{(n)}}\prec W^V_G$. It follows that $\psi(x,y,a)$ defines a relation without terminal nodes over $W^{\unionmodel{V}{n}}_{G^{(n)}}$. Let $\psi^*(x,y,a)$ be a relation defined over $\unionmodel{V}{n+1}$ by
\[
\textit{if } x\in W^{\unionmodel{V}{n}}_{G^{(n)}}, \textit{ then } y\in W^{\unionmodel{V}{n}}_{G^{(n)}} \textit{ and } W^{\unionmodel{V}{n}}_{G^{(n)}}\models\psi(x,y,a).
\]
Clearly, $\psi^*$ is a relation without terminal nodes over $\unionmodel{V}{n+1}$. Thus, by the $\DC_\omega$-scheme in $\unionmodel{V}{n+1}$, there is in $\unionmodel{V}{n+1}$ a function $f$ on $\omega$ that is a sequence of $\omega$-many dependent choices over $\psi^*(x,y,a)$. But clearly, since every initial segment of $f$ is in $W^{\unionmodel{V}{n}}_{G^{(n)}}$, as it is closed under finite sequences, we have that for all $m<\omega$, $W^{\unionmodel{V}{n}}_{G^{(n)}}\models \psi(f\restrict m,f(m),a)$. Thus, by elementarity, $W^V_G\models \psi(f\restrict m,f(m),a)$, for every $m<\omega$, as well.
\end{proof}

In order to prove the next theorem, we need the following result of Blass, which appears as Theorem 3.6 in \cite{blass:countablyManyGenericReals}.
\begin{theorem}\label{th:blassCohenReals}
A forcing extension $V[G]$ by $\Add(\omega,1)$ cannot have a sequence $\la r_\alpha\mid\alpha<\omega_1\ra$ of Cohen reals such that for every $\alpha<\omega_1$, $r_\alpha$ is $V[\la r_\xi\mid\xi<\alpha\ra]$-generic for $\Add(\omega,1)$.
\end{theorem}
The proof we give here is a slight modification of Blass's proof that will allow us to generalize the result in the next section.
\begin{proof}
Let $\mathbb B$ be the Boolean completion of $\Add(\omega,1)$. In particular, $\mathbb B$ has a dense subset of size $\omega$. Now suppose towards a contradiction that a forcing extension by $\mathbb B$ (equivalentely $\Add(\omega,1))$ has a sequence $\la r_\alpha\mid \alpha<\omega_1\ra$ of Cohen reals such that for every $\alpha<\omega_1$, $r_\alpha$ is Cohen generic over $V[\la r_\xi\mid\xi<\alpha\ra]$. The model $V[\la r_\xi\mid\xi<\omega_1\ra]$ is a forcing extension of $V$ by a complete subalgebra, $\mathbb D$, of $\mathbb B$ by the Intermediate Model Theorem of Solovay (see \cite{Grigorieff:intermediatemodels}). 

Let's first argue that $\mathbb D$ also has a dense subset of size $\omega$. Given a condition $p\in \Add(\omega,1)$, let $q_p$ be the infima of $b$ in $\mathbb D$ such that $p\leq b$. Each $q_p$ is in $\mathbb D$ by completeness and the conditions $q_p$ are dense in $\mathbb D$.  

Next, let $\dot R$ be a $\mathbb D$-name such that it is forced by $\one$ that $\dot R$ is an $\omega_1$-sequence of successively more generic Cohen reals and the extension by $\mathbb D$ is equal to the extension $V[\dot R]$. We claim that the Boolean values $\llbracket n\in \dot R(\alpha)\rrbracket$ for $n<\omega$ and $\alpha<\omega_1$ must generate $\mathbb D$. Suppose to the contrary that they generates a proper subalgebra $\mathbb D'$ of $\mathbb D$. Let $G$ be any $V$-generic filter for $\mathbb D$. Since $n\in \dot R_G(\alpha)$ if and only if $\llbracket n\in \dot R(\alpha)\rrbracket\in G$, we have that $\dot R_G$ already exists in $V[G\cap\mathbb D']$, which is a proper submodel of $V[G]=V[\dot R_G]$. 

Finally, observe that since $\mathbb D$ has a countable dense subset, there must be some $\alpha<\omega_1$ such that $\mathbb D$ is generated by the Boolean values $\llbracket n\in \dot R(\xi)\rrbracket$ for $n<\omega$ and $\xi<\alpha$. But this means that if $V[G]$ is $\mathbb D$-generic, then $G$ can be recovered from the sequence $\la r_\xi\mid\xi<\alpha\ra$, which contradicts that $r_{\alpha+1}$ is $V[\la r_\xi\mid\xi<\alpha\ra]$-generic.
\end{proof}

\begin{theorem}\label{th:DC_omegaHoldsDC_omega_2Fails}
Suppose that $V\models{\rm CH}$ and $\p=\Add(\omega,1)$ is the Cohen poset. Then
\begin{enumerate} \setlength \itemsep{2pt}
\item $W^V_G\models\ZFCminus+\DC_\omega$-scheme.
\item $W^V_G$ has the same cardinals and cofinalities as $V$.
\item $\mathcal P (\omega)$ is not $\omega_2$-big in $W^V_G$.
\item The $\DC_{\omega_2}$-scheme fails in $W^V_G$.
\item If additionally $V = L$, then the $\DC_{\omega_1}$-scheme fails in $W^V_G$.
\end{enumerate}
\end{theorem}

\begin{proof}
Clearly, $\Add(\omega,1)\cong \suppprod{\Add(\omega,1)}{\omega}{\omega}$. Item (1) follows from the theorems of Zarach above. Since $\Add(\omega,1)$ has the ccc, $V$ and $V[G]$ have the same cardinals and cofinalities, and hence so do $V$ and $W^V_G\subseteq V[G]$. Clearly, by ${\rm CH}$ in $V[G]$, $W^V_G$ cannot have a surjection from $\mathcal P (\omega)$, a proper class in $W^V_G$, onto $\omega_2$. Thus, by  \myref{prop:bigClasses}, the $\DC_{\omega_2}$-scheme fails in $W^V_G$. 

Now assume that $V = L$. The crucial observation is that this implies that $V$ is definable in $W^V_G$. In $W^V_G$, let $\varphi(x,y)$ be a formula asserting that whenever $x$ is a sequence of $L$-generic Cohen reals of length some $\alpha<\omega_1$, then $y$ is $L[x]$-generic for $\Add(\omega,1)$. The relation defined by $\varphi(x,y)$ has no terminal nodes because any sequence $x$ of $L$-generic Cohen reals is an element of some $V[G_n]$, and so $y$ given by $G_{\{n+1\}}$ works. Thus, if the $\DC_{\omega_1}$-scheme held in $W^V_G$, we would get an $\omega_1$-sequence of $L$-generic Cohen reals, which would contradict \myref{th:blassCohenReals}.
\end{proof}

\section{Generalized union models}\label{sec:GeneralizedZarach}

\noindent In this section, we will generalize Zarach's construction using products
\[
\q=\suppprod{\p}{\mu}{\delta}
\]
for regular cardinals $\delta\leq \mu$ to obtain failures of the $\DC_\delta$-scheme for larger cardinals $\delta$. The construction generalizes in a straightforward manner so we will just summarize the results here.

Suppose that $V\models\ZFC$. Let $\p\in V$ be a poset and let $\delta\leq\mu$ be regular cardinals such that $\q=\suppprod{\p}{\mu}{\delta}$ is isomorphic to $\p$. Let us call an automorphism $\pi$ of $\q$ \emph{coordinate-switching} if, as before, it acts by switching coordinates according to some automorphism $\treeauto{\pi}$ of $\mu$.

Let $G\subseteq \q$ be $V$-generic. For $\xi<\mu$, let
\begin{enumerate}
\item $G_\xi$ be the restriction of $G$ to the first $\xi$ coordinates,
\item $G_{\{\xi\}}$ be the restriction of $G$ to the $\xi$-th coordinate,
\item $G_{\xi,\tail}$ be the tail of $G$ after $\xi$ ($G=G_\xi\times G_{\xi,\tail}$).
\end{enumerate}
Let $\unionmodel{V}{\xi}=V[G_\xi]$. Let $W^V_G=\bigcup_{\xi<\mu}\unionmodel{V}{\xi}$.

Let $W^V_G$ be defined in $V[G]$ by the formula $\varphi(x,\la\p,\mu,\delta\ra,G)$, and for every formula $\psi(x)$ and $a\in W^V_G$, let the formula $\psi^W(x,\la \p,\mu,\delta\ra,G)$ be such that $W^V_G\models\psi(x)$ if and only if $V[G]\models \psi^W(a,\la \p,\mu,\delta\ra,G)$. As before, if $\pi$ is a coordinate-switching automorphism of $\q$, then both $\varphi(x,\la \p,\mu,\delta\ra,G)$ and $\varphi(x,\la \p,\mu,\delta\ra,\pi\image G)$ define $W^V_G$ in $V[G]$. Also, for every formula $\psi(x)$, we have that for every $a\in W^V_G$, $V[G]$ satisfies $\psi^W(a,\la \p,\mu,\delta\ra,G)$ if and only if $V[G]$ satisfies $\psi^W(a,\la \p,\mu,\delta\ra,\pi\image G)$. By an analogous automorphism argument as before, we get that if some condition $p\forces\psi^W(\check a,\la \check{\p},\check{\mu},\check{\delta}\ra,\dot G)$, where $\dot G$ is the canonical $\q$-name for the generic filter, then $\one\forces\psi^W(\check a,\la \check{\p},\check\mu,\check\delta\ra,\dot G)$.

As before, we shall write $W^V_G$ as the union of a sequence of models $W^{\unionmodel{V}{\xi}}_{G^{(\xi)}}$, each of which is an elementary submodel of $W^V_G$. To do this, fix an isomorphism $h \colon \p\cong \q=\suppprod{\p}{\mu}{\delta}$. Define the $V$-generic filters $G^{(\xi)}=h\image G_{\{\xi\}}$, for $\xi<\mu$. From $G^{(\xi)}$, for $\nu<\mu$, we obtain $G^{(\xi)}_{\nu}$, $G^{(\xi)}_{\{\nu\}}$, and $G^{(\xi)}_{\nu,\tail}$. Let
\[
\unionmodel{(\unionmodel{V}{\xi})}{\nu}=V[G_\xi][G^{(\xi)}_\nu]=V_{[\xi]}[G^{(\xi)}_\nu]
\]
for $\nu<\mu$, and let
\[
W^{\unionmodel{V}{\xi}}_{G^{(\xi)}}=\bigcup_{\nu<\mu} \unionmodel{(\unionmodel{V}{\xi})}{\nu}.
\]
Thus, we have:
\[
\unionmodel{V}{\xi} \subseteq \unionmodel{(\unionmodel{V}{\xi})}{1} \subseteq \unionmodel{(\unionmodel{V}{\xi})}{2} \subseteq \cdots \subseteq \unionmodel{(\unionmodel{V}{\xi})}{\nu} \subseteq \cdots \subseteq \bigcup_{\nu<\mu} \unionmodel{(\unionmodel{V}{\xi})}{\nu} = W^{\unionmodel{V}{\xi}}_{G^{(\xi)}} \subseteq \unionmodel{V}{\xi + 1}.
\]

\noindent An analogous argument to the proof of \myref{lem:elementarityW} yields.

\begin{lemma}\label{lem:elementarityWGeneral}
For every $\xi<\mu$, $W^{\unionmodel{V}{\xi}}_{G^{(\xi)}}\prec W^V_G$.
\end{lemma}

\noindent Observe that $W^V_G$ cannot be a model of the power set axiom because $\mathcal P (\p)$ does not exist in $W^V_G$. \myref{lem:elementarityWGeneral} gives:

\begin{theorem}
$W^V_G\models\ZFCminus$.
\end{theorem}

\begin{theorem}\label{th:closureOfW}
Suppose that $\p$ is ${<}\delta$-closed. Then $(W^V_G)^{{<}\delta}\subseteq W^V_G$ in $V[G]$.
\end{theorem}

\begin{proof}
Since $\p$ is ${<}\delta$-closed, the product $\q$ is ${<}\delta$-closed as well. Since $\mu\geq\delta$ was regular in $V$, it follows that the cofinality of $\mu$ in $V[G]$ is at least $\delta$. Suppose that $f \colon \gamma \to W^V_G$ in $V[G]$ for some $\gamma<\delta$. Then there is some $\xi<\mu$ such that the range of $f$ is contained in $\unionmodel{V}{\xi}=V[G_\xi]$ by cofinality considerations. But since the tail of the product after $\xi$ is ${<}\delta$-closed, it follows that $f\in \unionmodel{V}{\xi}$, and hence $f\in W^V_G$.
\end{proof}

\begin{theorem}
Suppose that $\p$ is ${<}\delta$-closed. Then $W^V_G\models\DC_\delta$-scheme.
\end{theorem}

\begin{proof}
Fix a formula $\psi(x,y,a)$ defining over $W^V_G$ a relation without terminal nodes. Let $\xi$ be large enough so that $a\in \unionmodel{V}{\xi}$. By \myref{lem:elementarityWGeneral}, $W^{\unionmodel{V}{\xi}}_{G^{(\xi)}}\prec W^V_G$. It follows that $\psi(x,y,a)$ defines a relation without terminal nodes over $W^{\unionmodel{V}{\xi}}_{G^{(\xi)}}$. Let $\psi^*(x,y,a)$ be a relation defined over $\unionmodel{V}{\xi + 1}$ by whenever $x\in W^{\unionmodel{V}{\xi}}_{G^{(\xi)}}$, then $y\in W^{\unionmodel{V}{\xi}}_{G^{(\xi)}}$ and $W^{\unionmodel{V}{\xi}}_{G^{(\xi)}}\models\psi(x,y,a)$. Thus, by the $\DC_\delta$-scheme in $\unionmodel{V}{\xi + 1}$, there is in $\unionmodel{V}{\xi + 1}$ a function $f$ on $\delta$ that is a sequence of $\delta$-many dependent choices over $\psi^*(x,y,a)$. Now we use the ${<}\delta$-closure of $W^{\unionmodel{V}{\xi}}_{G^{(\xi)}}$ (\myref{th:closureOfW}) to confirm that every initial segment of $f$ is in $W^{\unionmodel{V}{\xi}}_{G^{(\xi)}}$. Thus, for all $\nu<\delta$, $W^{\unionmodel{V}{\xi}}_{G^{(\xi)}}\models \psi(f\restrict \nu,f(\nu),a)$. So now by elementarity, $W^V_G\models \psi(f\restrict \nu,f(\nu),a)$, for every $\nu<\delta$, as well.
\end{proof}

\noindent Given a regular cardinal $\delta$, let $\Add(\delta,1)$ be the generalized Cohen poset adding a subset to $\delta$ with conditions of size less than $\delta$. First, we state a generalization of \myref{th:blassCohenReals} from the previous section. 

\begin{theorem}\label{th:blassCohenSubsets}
Suppose $\delta$ is a regular cardinal with $\delta^{{<}\delta}=\delta$. A forcing extension $V[G]$ by $\Add(\delta,1)$ cannot have a sequence $\la A_\alpha\mid\alpha<\delta^+\ra$ of Cohen subsets such that for every $\alpha<\delta^+$, $A_\alpha$ is $V[\la A_\xi\mid\xi<\alpha\ra]$-generic for $\Add(\delta,1)$.
\end{theorem}

The proof is completely analogous to the proof of \myref{th:blassCohenReals}, using the assumption $\delta^{{<}\delta}=\delta$ to show that $\Add(\delta,1)$ has size $\delta$.

\pagebreak
\begin{theorem}\label{th:DC_deltaHOldsDC_delta++Fails}
Suppose that $V\models2^\delta=\delta^+$ for some regular cardinal $\delta$ and let $\p=\Add(\delta,1)$. Then
\begin{enumerate} \setlength \itemsep{1pt}
\item $W^V_G\models \ZFCminus + \DC_\delta$-scheme.
\item $W^V_G$ has the same cardinals and cofinalities as $V$, with the possible exception of $\delta^+$.
\item $\mathcal P (\delta)$ is not $\delta^{++}$-big in $W^V_G$.
\item The $\DC_{\delta^{++}}$-scheme fails in $W^V_G$.
\item If additionally $V = L$, then the $\DC_{\delta^+}$-scheme fails in $W^V_G$.
\end{enumerate}
\end{theorem}

\begin{proof}
Clearly, $\Add(\delta,1)\cong \suppprod{\Add(\delta,1)}{\delta}{\delta}$. Item (1) follows from the theorems above. Since $\Add(\delta,1)$ is ${<}\delta$-closed and has at most $\delta^{++}$-cc (by $2^\delta=\delta^+$), $V$ and $V[G]$ have the same cardinals and cofinalities with the possible exception of $\delta^+$, and hence so do $V$ and $W^V_G\subseteq V[G]$. Clearly, since $2^\delta=\delta^+$ holds in $V[G]$, $W^V_G$ cannot have a surjection from $\mathcal P (\delta)$, a proper class in $W^V_G$, onto $\delta^{++}$. Thus, by \myref{prop:bigClasses}, the $\DC_{\delta^{++}}$-scheme fails in $W^V_G$. If $V = L$, then $\delta^{{<}\delta}=\delta$ and the $\DC_{\delta^+}$-scheme fails by an application of \myref{th:blassCohenSubsets} as in the proof of \myref{th:DC_omegaHoldsDC_omega_2Fails}.
\end{proof}

\section{A large model where the \texorpdfstring{$\DC_\omega$-scheme}{DC omega scheme} fails}\label{sec:LargeModelNegDC}

\noindent In this section we shall provide a union model in the style of Zarach for which the $\DC_\omega$-scheme fails. Unlike the small model of \cite{friedmangitmankanovei19}, this model will have unboundedly many cardinals.

We work in the second-order model $\mathcal V=\la L,\in,\mathcal C\ra$, where $\C$ is the collection of definable classes of $L$. We will force with the class tree iteration $\p(\vec \J,\Ord^{{<}\omega})$. Let\break $G\subseteq \p(\vec \J,\Ord^{{<}\omega})$ be $\mathcal V$-generic. By \myref{prop:JensenTreeForcingPretame}, $\p(\vec \J,\Ord^{{<}\omega})$ is pretame, and hence $\mathcal V[G]=\la L[G],\in,\mathcal C[
G]\ra\models \GBcminus$. Although, we won't make use of this fact, let's also note that $\mathcal V[G]\models\DCOrd$-scheme by \myref{prop:pretamePreservesDC}.

Extending our earlier terminology, we will call an automorphism $\pi$ of $\p(\vec \J,\Ord^{{<}\omega})$ \emph{tree-switching} if there is an automorphism $\treeauto{\pi}$ of $\Ord^{{<}\omega}$ such that for any condition $p$ $\pi(p)=q$, where $q(t)=p(\treeauto{\pi}^{-1}(t))$, namely $\pi$ switches the nodes of $\Ord^{{<}\omega}$ according to $\treeauto{\pi}$.

Fix a set tree $\mathcal T\subseteq\Ord^{{<}\omega}$. Let $G_{\mathcal T}\subseteq G$ consist of all functions $f_T\in \p(\vec \J,\Ord^{{<}\omega})$ with $T\subseteq \mathcal T$ a finite subtree. Let's argue that $G_{\mathcal T}$ is $L$-generic for $\p(\vec \J,\mathcal T)$. It suffices to show that every maximal antichain $A$ of $\p(\vec \J,\mathcal T)$ remains maximal in $\p(\vec \J, \Ord^{{<}\omega})$. Fix a maximal antichain $A$ of $\p(\vec \J,\mathcal T)$. Take any $f_S\in \p(\vec \J,\Ord^{{<}\omega})$. Let $S^*=S\cap \mathcal T$ and let $f_{S^*}=f_S\restrict S^*$, so that $f_{S^*}\in \p(\vec \J,\mathcal T)$. By the maximality of $A$ in $\p(\vec \J,\mathcal T)$, there is $f_T\in A$ compatible with $f_{S^*}$. But then clearly $f_T$ is compatible with $f_S$ as well. Thus, $G_\mathcal{T}$ is an $L$-generic for $\p(\vec \J,\mathcal T)$ and therefore $L[G_{\mathcal T}]\models\ZFC$.

Let $\mathbb T$ consist of all infinite trees $\mathcal T\subseteq \Ord^{{<}\omega}$ such that $\mathcal T$ does not have a cofinal branch. Let $W=\bigcup_{\mathcal T\in \mathbb T}L[G_{\mathcal T}]$. We will show below that $W\models\ZFCminus+\neg\DC_{\omega}$-scheme. But first we need some technical preliminaries.

\begin{proposition}\label{prop:isomorphicSubtreesExtensions}
Suppose that $\mathcal S_1$ and $\mathcal S_2$ are subtrees of $\Ord^{{<}\omega}$, and $\pi$ is a \break tree-switching automorphism of $\p(\vec \J,\Ord^{{<}\omega})$ such that $\treeauto{\pi} \image \mathcal S_1 = \mathcal S_2$. Then \break $L[(\pi^{-1}\image G)_{\mathcal S_1}]=L[G_{\mathcal S_2}]$.
\end{proposition}

\begin{proof}
The automorphism $\pi^{-1}$ restricts to an isomorphism from $\p(\vec \J,\mathcal S_2)$ to $\p(\vec \J,\mathcal S_1)$ and the image of $G_{\mathcal S_2}$ under this isomorphism is $(\pi^{-1}\image G)_{\mathcal S_1}$. Thus,
\[
L[(\pi^{-1}\image G)_{\mathcal S_1}] \cong L[G_{\mathcal S_2}],
\]
and hence $L[(\pi^{-1}\image G)_{\mathcal S_1}]=L[G_{\mathcal S_2}]$.
\end{proof}

\noindent Observe that, since $L$ is definable in $L[G]$, there is a formula $\varphi(x,\vec \J,G)$ defining $W$ in $L[G]$. Thus, for every formula $\psi(x)$, there is a corresponding formula $\psi^W(x,y,z)$ such that for every $a\in W$, $W\models\psi(a)$ if and only if $L[G]\models\psi^W(a,\vec\J,G)$. Next, let's argue that if $\pi$ is any tree-switching automorphism of $\p(\vec \J,\Ord^{{<}\omega})$, then $\varphi(x,\vec\J,\pi\image G)$ also defines $W$ in $L[G]$. Fix $\mathcal T\in \mathbb T$. By \myref{prop:isomorphicSubtreesExtensions}, $L[(\pi\image G)_{\treeauto{\pi}^{-1} \image \mathcal T}]=L[G_{\mathcal T}]$. Thus, $\bigcup_{\mathcal T\in\mathbb T}L[G_{\mathcal T}]=\bigcup_{\mathcal T\in\mathbb T}L[(\pi\image G)_{\mathcal T}]$. Hence also, for every formula $\psi(x)$, we have that for every $a\in W$, $W\models\psi(a)$ if and only if $L[G]\models\psi^W(a,\vec \J,\pi\image G)$.

\begin{proposition}\label{prop:RestrictedConditionForces}
Suppose that for some formula $\psi(x)$, $p\forces\psi^W(\dot a)$, where $\dot a$ is a $\p(\vec\J,\mathcal T)$-name, for some tree $\mathcal T\in \mathbb T$. Then $p\restrict \mathcal T\forces\psi^W(\dot a)$.
\end{proposition}

\begin{proof}
Suppose towards a contradiction that $p\restrict\mathcal T$ does not force $\psi^W(\dot a)$. Then there is a condition $q\leq p \restrict \mathcal T$ such that $q\forces \neg\psi^W(\dot a)$. Let $\pi$ be a tree-switching automorphism such that $\treeauto{\pi}$ fixes $\mathcal T$ and moves the nodes in $\text{dom}(p) \setminus \mathcal T$ so that
\[
\text{dom}(q)\cap \text{dom}(\pi(p)) = \text{dom}(p \restrict \mathcal T) \subseteq \mathcal T.
\]
We have $\pi(p)\forces\psi^W(\dot a)$ because $\pi(\dot a)=\dot a$ and tree-switching automorphisms don't affect $W$. But this is impossible because clearly $\pi(p)$ and $q$ are compatible.
\end{proof}

\begin{theorem}\label{th:WJensenZFC-}
$W\models \ZFCminus$.
\end{theorem}

\begin{proof}
It is clear that $W$ satisfies extensionality, empty set, pairing, unions, infinity, the foundation scheme, and the well-ordering principle. We will be done if we can argue that $W$ satisfies the replacement and collection schemes (separation will then follow). We will verify collection because the same argument will yield replacement as well.

Since $W$ satisfies the well-ordering principle, it suffices to verify instances of collection for ordinals. So suppose that
\[
W\models\forall \xi<\delta\exists y\,\psi(\xi,y,a).
\]
Let $a\in L[G_{\mathcal T}]$ for some $\mathcal T\in\mathbb T$, and let $\dot a$ be a $\p(\vec\J,\mathcal T)$-name for $a$. Let \hbox{$p \forces \theta^W(\dot{a})$,} where $\theta(\dot{a}) := \forall \xi< \check\delta \exists y\,\psi(\xi,y,\dot{a})$. By \myref{prop:RestrictedConditionForces}, we can assume without loss of generality that $\dom(p)\subseteq\mathcal T$. Given a tree $\mathcal S\in \mathbb T$, let $\dot G_{\mathcal S}$ be the canonical $\p(\vec \J,\mathcal S)$-name for the generic filter. Observe that if $\pi$ is a tree-switching automorphism, then $\pi(\dot G_{\mathcal S})=\dot G_{\treeauto{\pi} \image \mathcal S}$.

Before giving the details of the proof, we sketch the idea behind the argument: For each $\xi < \delta$ there is some tree $\mathcal S_\xi \in \mathbb T$ such that $W \models \exists y \in L[G_{\mathcal S_\xi}] \psi(\xi, y, a)$. The aim is to find some tree $\mathcal R \in \mathbb T$ such that, for each $\xi \in \delta$, $W \models \exists y \in L[G_{\mathcal R}] \psi(\xi, y, a)$.

Since $\bigcup_{\xi < \delta} S_\xi$ may contain a cofinal branch, we want to find a sequence of tree-switching automorphisms $\pi_\xi$ such that $W \models \exists y \in L[G_{\pi_\xi \image \mathcal S_\xi}] \psi(\xi, y, a)$ and $\pi_\xi \image \mathcal S_\xi$ and $\pi_\eta \image \mathcal S_\eta$ are disjoint modulo $\mathcal T$ for any $\xi, \eta < \delta$.

However, we are unable to determine these automorphisms in the ground model and therefore it need not be the case that $\bigcup_{\xi < \delta} \pi_\xi \image \mathcal S_\xi \in \mathbb{T}$ (in $L$). To avoid this issue we shall use the ccc to recursively construct in $L$ countable sequences of trees $\langle S_\xi^{(\alpha)} \mid \alpha < \beta_\xi \rangle$ for $\xi < \delta$ such that $S^{(\alpha)}_\xi$ and $S^{(\gamma)}_\eta$ are pairwise disjoint module $\mathcal T$ for every $\alpha, \xi, \gamma$ and $\eta$. This will be done in such a way that, for $\mathcal T_\xi = \bigcup_{\alpha < \beta_\xi} S_\xi^{(\alpha)}$,
\[
W \models \forall \xi < \delta \exists y \in L[G_{\mathcal T_\xi}] \psi(\xi, y, a).
\]
We will then be able to take $\mathcal R := \bigcup_{\xi \in \delta} \mathcal T_\xi \in \mathbb T$ as our witnessing tree. So that this can be easily modified into a proof for replacement, we will also explicitly construct names $y_\xi$ for each $\xi < \delta$.

\bigskip

\noindent In order to do this, for every $\xi<\delta$, let $D_\xi$ be the dense class of conditions $q$ below $p$ for which there is some $\mathcal S \in\mathbb T$ such that $\mathcal T \subseteq \mathcal S$ and a $\p(\vec\J,\mathcal S)$-name $\dot y$ such that
\[
q\forces \dot y\in L[\dot G_{\mathcal S}]\wedge \psi^W(\check \xi,\dot y,\dot a).
\]
Choose any condition $q_0^{(0)}\in D_0$ (using global choice in $L$) and fix $\dot{y}^{(0)}_0$ and $\mathcal S^{(0)}_0$ such that
\[
q^{(0)}_0\forces \dot y_0^{(0)} \in L[\dot G_{\mathcal S_0^{(0)}}]\wedge \psi^W(\check 0,\dot y^{(0)}_0,\dot a).
\]

\noindent Next, assuming this is possible, choose any condition $q_1^*\in D_0$ incompatible with $q_0^{(0)}$, and fix $\dot{y}^{*}_1$ and $\mathcal S^{*}_1$ such that
\[
q_1^*\forces \dot y_1^* \in L[\dot G_{\mathcal S_1^*}]\wedge\psi^W(\check 0,\dot y_1^*,\dot a).
\]

\noindent Let $\pi_1$ be a tree-switching automorphism such that $\treeauto{\pi}_1$ fixes $\mathcal T$ and $\mathcal S^{(1)}_0 := \treeauto{\pi}_1 \image \mathcal S^{*}_1$ satisfies $\mathcal S_0^{(1)} \cap \mathcal S_0^{(0)} \subseteq \mathcal T$. Now, we have that
\[
\pi_1(q_1^*)\forces \pi_1(\dot y_1^*)\in L[\dot G_{\mathcal S_0^{(1)}}]\wedge \psi^W(\check 0,\pi_1(\dot y_1^*),\dot a).
\]
Also, since $\dom(p)\subseteq \mathcal T$ and $\treeauto{\pi}_1$ fixes $\mathcal T$, $\pi_1(q_1^*)\leq p$. Let $q_0^{(1)}=\pi_1(q_1^*)\cup q_1^*$. Thus, it is still case that $q_0^{(1)}\in D_0$ and $q_0^{(1)}$ is not compatible with $q_0^{(0)}$. Let $\dot y_0^{(1)}=\pi_1(\dot y_1^*)$. Continuing in this manner, we keep building a sequence of mutually incompatible conditions $q_0^{(\alpha)}\in D_0$ such that
\[
q_0^{(\alpha)}\forces \dot y_0^{(\alpha)}\in L[\dot G_{\mathcal S_0^{(\alpha)}}]\wedge \psi^W(\check 0,\dot y_0^{(\alpha)},\dot a)
\]
and $\mathcal S_0^{(\alpha)} \cap \bigcup_{\gamma < \alpha} \mathcal S_0^{(\gamma)} \subseteq \mathcal T$. This process must terminate after $\beta_0$-many steps for a countable $\beta_0$ because the poset $\p(\vec \J,\Ord^{{<}\omega})$ has the ccc. Let $A_0=\{q_0^{(\alpha)} \mid \alpha < \beta_0\}$ be the resulting maximal antichain contained in $D_0$. Let $\mathcal T_0 = \bigcup_{\alpha < \beta_0}\mathcal S_0^{(\alpha)}$, and observe that by the disjointness of the $S_0^{(\alpha)}$ modulo $\mathcal T$, we have that $\mathcal T_0$ cannot have an infinite branch and therefore $\mathcal T_0 \in \mathbb T$.

Let $\dot z_0$ be the mixed name of the names $\dot y_0^{(\alpha)}$ over the antichain $A_0$. Namely,
\[
\dot z_0 = \bigcup_{\alpha < \beta_0} \Big\{ \la \dot x, r \ra \mid r \leq q_0^{(\alpha)}, r \in \p(\vec\J, \mathcal S_0^{(\alpha)}), r \forces \dot x \in \dot y_0^{(\alpha)}, \dot x \in \text{dom}(\dot y_0^{(\alpha)}) \Big\}.
\]

\noindent Note that we can include the condition $r \in \p(\vec\J, \mathcal S_0^{(\alpha)})$ because if $r \forces \dot x \in \dot y_0^{(\alpha)}$ for some $\dot x \in \text{dom}(\dot y_0^{(\alpha)})$, then $r \restrict \mathcal S_0^{(\alpha)}$ also forces this by \myref{prop:RestrictedConditionForces}. Finally, observe that $\dot z_0$ is a $\p(\vec\J,\mathcal T_0)$-name and $W \models (\dot{z}_0)_G \in L[G_{\mathcal T_0}] \wedge \psi(\check 0, (\dot{z}_0)_G, a)$.

\bigskip

\noindent Next, we repeat the process for $D_1$, building a maximal antichain
\[
A_1=\{q_1^{(\alpha)} \mid \alpha < \beta_1\}
\]
contained in $D_1$ and trees $\mathcal S_1^{(\alpha)}$ such that
\[
q_1^{(\alpha)} \forces \dot y_1^{(\alpha)} \in L[\dot G_{\mathcal S_1^{(\alpha)}}] \wedge \psi^W(\check 1,\dot y_1^{(\alpha)},\dot a).
\]
At the same time, we ensure that for any $\alpha < \beta_1$, $\mathcal S_1^{(\alpha)} \cap (\bigcup_{\gamma < \alpha} \mathcal S_1^{(\gamma)} \cup \mathcal T_0) \subseteq \mathcal T$. Let $\mathcal T_1 = \bigcup_{\alpha < \beta_1} \mathcal S_1^{(\alpha)}$ and observe that $\mathcal T_0 \cup \mathcal T_1$ is in $\mathbb T$. Let $\dot z_1$ be the mixed name of the names $\dot y_1^{(\alpha)}$ over the antichain $A_1$, and observe that $\dot z_1$ is a $\p(\vec \J,\mathcal T_1)$-name.

We continue the process for every $D_\xi$, and let $\mathcal R = \bigcup_{\xi < \delta} \mathcal T_\xi$, which is in $\mathbb T$ by construction. Let $\dot z$ be the canonical name for a sequence of length $\delta$ obtained from the names $\dot z_\xi$ for $\xi<\delta$. Then $\dot z$ is a $\p(\vec\J,\mathcal R)$-name. By construction, for every $\xi<\delta$, $W\models \psi(\check \xi,\dot z_G(\xi),a)$, so $\dot z_G$ witnesses this instance of collection.
\end{proof}

\begin{theorem}\label{th:WJensennotReflection}
$W\models\neg \DC_\omega$-scheme.
\end{theorem}

\begin{proof}
Consider the definable class tree whose domain is
\[
\{\vec r \mid \vec r \text{ is $L$-generic for }\J_n\text{ for some $n$}\}
\]
ordered by extension in $W$. Clearly, the tree relation has no terminal nodes. Thus, if we can show that it doesn't have an infinite branch, we will have a violation of the $\DC_\omega$-scheme. So suppose that $b\in W$ is an infinite branch through this class tree. Then $b\in L[G_{\mathcal S}]$ for some tree $\mathcal S\in \mathbb T$. Since $\mathcal S$ does not have an infinite branch by the definition of $\mathbb T$, there must be some $\vec r_n$, the element of $b$ on level $n$, which is $L$-generic for $\J_n$ but not in $\mathcal S$. However, this is impossible by \myref{th:JensenTreeForcing}~(2).
\end{proof}

\noindent The model $W$ is also interesting because even though the $\DC_{\omega}$-scheme fails, every proper class in $W$ is big. Before we prove this we need the following lemma.

\begin{lemma}\label{lem:intersectionOfTreeModels}
Suppose that $\mathcal S_1,\mathcal S_2,\mathcal T\in\mathbb T$ are such that $\mathcal S_1\cap \mathcal S_2=\mathcal T$. Then
\[
L[G_{\mathcal S_1}]\cap L[G_{\mathcal S_2}]=L[G_{\mathcal T}].
\]
\end{lemma}

\begin{proof}
Since we are dealing with models of $\ZFC$, it suffices to show that every set of ordinals in $L[G_{\mathcal S_1}]\cap L[G_{\mathcal S_2}]$ is in $L[G_{\mathcal T}]$. Suppose that $A\in L[G_{\mathcal S_1}]\cap L[G_{\mathcal S_2}]$ is a subset of an ordinal $\alpha$. Let $\dot x$ be a nice $\p(\vec \J,\mathcal S_1)$-name for $A$, namely $\dot x = \bigcup_{\xi<\alpha} \{\check \xi\} \times A_\xi$, where the $A_{\xi}$ are antichains of $\p(\vec\J,\mathcal S_1)$. Similarly, let $\dot y = \bigcup_{\xi<\alpha}\{\check \xi\} \times B_\xi$ be a nice $\p(\vec \J,\mathcal S_2)$-name for $A$. Fix a condition $p\in \p(\vec \J,\mathcal S_1\cup\mathcal S_2)$ forcing that $\dot x=\dot y$. Let $p=p_1\cup p_2$, where $p_1=p\restrict \mathcal S_1$ and $p_2=p\restrict \mathcal S_2$. We will work below this condition $p$. By shrinking the set $A$, we can assume without loss of generality that $p$ does not decide $\check \xi \in \dot x$ for any $\xi<\alpha$. We can also assume that conditions in all $A_\xi$ are compatible with $p$.

Next, let's argue that if some condition $q\leq p_1$ in $\p(\vec\J,\mathcal S_1)$ decides $\check \xi\in\dot x$, then $p_1\cup q\restrict \mathcal T$ already decides $\check\xi\in\dot x$. Suppose that $q \forces \check \xi\in \dot x$ (the case $q\forces \check\xi\not\in\dot x$ will be the same). We will first show that $p_2 \cup q \restrict \mathcal T$ forces that $\check \xi \in \dot y$. So, suppose that this is not the case. Then fix $r \leq p_2 \cup q \restrict \mathcal T$ in $\p(\vec\J,\mathcal S_2)$ such that $r \forces \check\xi \notin \dot y$. But then $q \cup r \leq q$ in $\p(\vec\J,\mathcal S_1\cup \mathcal S_2)$, and $q\cup r\forces \check\xi \notin \dot y$ in $\p(\vec\J,\mathcal S_1\cup \mathcal S_2)$ (by absoluteness for atomic forcing formulas), and also $q\cup r \forces \check\xi \in \dot x$ in $\p(\vec\J,\mathcal S_1 \cup \mathcal S_2)$. But this is a contradiction because $q \cup r \leq p_1 \cup p_2 \cup q \restrict \mathcal T \leq p$ and so must force $\dot x=\dot y$. Thus, $p_2 \cup q \restrict \mathcal T \forces \check\xi \in\dot y$. But now essentially the same argument on the $\mathcal S_1$-side with $\dot x$ shows that $p_1 \cup q \restrict \mathcal T \forces \check\xi \in \dot x$.

Let $\dot x^*=\bigcup_{\xi<\alpha}\{\check \xi\}\times A^*_\xi$, where $A^*_\xi=\{q\restrict \mathcal T\mid q\in A_\xi\}$, and note that $\dot x^*$ is a $\p(\vec\J,\mathcal T)$-name. We claim that $\dot x^*_G=\dot x_G$. Suppose that $\xi \in \dot x_G$. Then there is $q \in G\cap A_\xi$. Thus, $q \restrict \mathcal T \in G \cap A_\xi^*$, and hence $\xi \in \dot x^*_G$. Next, suppose that $\xi \notin \dot x_G$. Then there is $q\leq p$ in $G$ such that $q\forces \check \xi\notin \dot x$. By the above argument, the condition $p_1 \cup q \restrict \mathcal T$ also forces $\check \xi \notin \dot x$. Thus, conditions incompatible with some $a\in A_\xi$ are dense below $p_1 \cup q \restrict \mathcal T$. But then, by our assumption that $p$, and therefore $p_1$, is compatible with all conditions in every $A_\xi$, it follows that conditions incompatible with some $a\in A_\xi^*$ are dense below $q\restrict\mathcal T$. Thus, $q\restrict \mathcal T\forces \check\xi \not\in \dot x^*$. Since $q\restrict\mathcal T\in G$, it follows that $\xi \notin \dot x^*_G$.
\end{proof}

\begin{theorem} \label{th:WEveryClassisBig}
Every proper class in $W$ is big.
\end{theorem}

\begin{proof}
Suppose that a formula $\psi(x,a)$ defines a proper class $\mathcal A$ in $W$ and $a\in L[G_{\mathcal T}]$ for some $\mathcal T\in \mathbb T$. Fix a cardinal $\delta$ and recall that, since $\p(\vec \J,\Ord^{{<}\omega})$ has the ccc, all of our models have the same cardinals. We need to verify that there is a surjection from $\mathcal A$ onto $\delta$. First, suppose that $\mathcal A \cap L[G_{\mathcal T}]$ is a proper class. In this case, there must be some ordinal $\alpha$ such that $\mathcal A\cap L[G_{\mathcal T}]_\alpha$ has size at least $\delta$ in $L[G]$. Let $A=\mathcal A\cap L[G_{\mathcal T}]_\alpha$, which exists in $W$ by separation. Since $W$ can enumerate $A$ by the well-ordering principle, let $f \colon \beta \to A$ for some ordinal $\beta$, and observe that by our assumption on the size of $A$, $\beta \geq \delta$.

So let's assume now that $\mathcal A \cap L[G_{\mathcal T}]$ is not a proper class. As in \myref{th:WJensenZFC-}, we begin by sketching the idea behind the argument. Since $\mathcal A$ is a proper class in $W$, there will be some tree $\mathcal S_0$ extending $\mathcal T$ for which \mbox{$W \models \exists x \in L[G_{\mathcal S_0}] \setminus L[G_{\mathcal T}] \, \psi(x, a)$.} The aim is to find some sequence $\langle \mathcal S_\xi \mid \xi \in \delta \rangle$ each of which is isomorphic to $\mathcal S_0$, pairwise disjoint modulo $\mathcal T$, and such that
\[
W \models \exists x \in L[G_{\mathcal S_\xi}] \setminus L[G_{\mathcal T}] \, \psi(x, a).
\]
However, since these trees cannot be determined in the ground model we will again use the ccc to recursively construct countable sequences of trees $\langle S_\xi^{(\alpha)} \mid \alpha < \beta_\xi \rangle$ for $\xi < \delta$ such that $S_\xi^{(\alpha)} \cap S_\eta^{(\gamma)} \subseteq \mathcal T$ for each $\alpha, \xi, \gamma$ and $\eta$ and, for $\mathcal T_\xi = \bigcup_{\alpha < \beta_\xi} S_\xi^{(\alpha)}$,
\[
W \models \exists x \in L[G_{\mathcal T_\xi}] \setminus L[G_{\mathcal T}] \, \psi(x, a).
\]

So, let $\dot a$ be a $\p(\vec \J,\mathcal T)$-name for $a$. Fix a condition
\[
p \forces \exists x \notin L[\dot G_{\mathcal T}]\,\psi^W(x,\dot a).
\]
By \myref{prop:RestrictedConditionForces}, we can assume without loss of generality that $p \in \p(\vec\J,\mathcal T)$. Let $D$ be the dense class of conditions $q$ below $p$ forcing for some tree $\mathcal S\in\mathbb T$ that $\exists x\in L[\dot G_{\mathcal S}]\,\psi^W(x,\dot a)$. Following the proof of \myref{th:WJensenZFC-}, build a maximal antichain $A_0 = \{q_0^{(\alpha)} \mid \alpha < \beta_0\}$ contained in $D$ such that
\[
q_0^{(\alpha)}\forces \exists x\in L[\dot G_{\mathcal S_0^{(\alpha)}}]\,\psi^W(x,\dot a),
\]
where $\mathcal S_0^{(\alpha)} \cap \bigcup_{\gamma<\alpha}\mathcal S_0^{(\gamma)}\subseteq \mathcal T$. Let $\mathcal T_0 = \bigcup_{\alpha < \beta_0} S_0^{(\alpha)}$. Next, we repeat the process, constructing a maximal antichain $A_1 = \{q_1^{(\alpha)}\mid \alpha < \beta_1\}$ contained in $D$ such that
\[
q_1^{(\alpha)} \forces \exists x\in L[\dot G_{\mathcal S_1^{(\alpha)}}]\,\psi^W(x,\dot a),
\]
and we have $\mathcal S_1^{(\alpha)}\cap (\bigcup_{\gamma < \alpha}\mathcal S_1^{(\gamma)}\cup \mathcal T_0) \subseteq \mathcal T$. We continue this process, constructing maximal antichains $A_\xi = \{q_\xi^{(\alpha)} \mid \alpha < \beta_\xi\}$ contained in $D$ such that
\[
q_\xi^{(\alpha)} \forces \exists x \in L[\dot G_{\mathcal S_\xi^{(\alpha)}}]\,\psi^W(x,\dot a),
\]
maintaining the disjointness of the trees modulo $\mathcal T$.

Let $\mathcal R = \bigcup_{\xi<\delta} \mathcal T_\xi$ which is in $\mathbb T$ by the the disjointness of the trees $S_\xi^{(\alpha)}$ modulo $\mathcal T$. Finally, by \myref{lem:intersectionOfTreeModels}, for each $\xi < \delta$
\[
A \cap L[G_{\mathcal T_\xi}] \setminus \bigcup_{\eta < \xi} A \cap L[G_{\mathcal T_\eta}] \neq \emptyset
\]
and therefore the model $L[G_{\mathcal R}]$ contains at least $\delta$-many elements of $\mathcal A$.
\end{proof}

\noindent We should point out that the construction given above fails if we replace $\Ord^{{<}\omega}$ with $\alpha^{{<}\omega}$ for some cardinal $\alpha$. Indeed, the model $W$ constructed analogously in a forcing extension $L[G]$ by an $L$-generic $G\subseteq \p(\vec \J,\alpha^{{<}\omega})$ fails to satisfy the following instance of collection. We analogously let $\mathbb T$ be the collection of all subtrees of $\alpha^{{<}\omega}$ of size $\alpha$ which do not have an infinite branch. In this case, $\mathbb T$ is a set from $L$, and hence in the union model $W$. The model $W$ satisfies that for every $\mathcal T\in\mathbb T$, there is a constellation of Jensen reals along $\mathcal T$. More formally, there is a map $F_{\mathcal T}$ with domain $\mathcal T$ such that nodes of length $n$ get mapped to sequences of $L$-generic reals for $\J_n$ and the sequences on longer nodes end-extend sequences on shorter nodes. Suppose towards a contradiction that there is a collecting set $C$ for this instance of collection. $C$ must then be in some $L[G_{\mathcal T}]$ with $\mathcal T \in \mathbb T$. Let $\mathcal S$ be a tree in $\mathbb T$ of rank higher than $\mathcal T$. This ensures that there is no tree isomorphism between $\mathcal T$ and a subtree of $\mathcal S$. But then by \myref{th:JensenTreeForcing}~(2), $L[G_{\mathcal T}]$ cannot contain the map $F_{\mathcal S}$.

We should also note that instead of forcing over $L$, we could have forced over $H_{\lambda}^L$ for some regular, uncountable $\lambda$. The model $H_{\lambda}^L \models \ZFCminus$ and for it $\p(\vec\J,\lambda^{{<}\omega})$ is a pretame forcing. Thus, if $G\subseteq \p(\vec\J,\lambda^{{<}\omega})$ is $H_{\lambda}$-generic, then $H_{\lambda}[G]\models\ZFCminus$. The rest of the arguments in this section then go through.

In \myref{th:WEveryClassisBig}, we proved that $W$ is a model of $\ZFCminus$ in which every proper class is big. While this is a very desirable property for our model to satisfy, we can make one final observation for this section. This is that, by combining the construction with that of \myref{sec:GeneralizedZarach}, it is possible to produce a model of $\ZFCminus$ in which the $\DC_\omega$-scheme fails and in which there are proper classes that are not big.

\begin{corollary}
It is possible to produce a model of $\ZFCminus$ with unboundedly many cardinals in which the $\DC_\omega$-scheme fails and there is a proper class that is not big.
\end{corollary}

\begin{proof}
We start with a model $\mathcal V$ of $\GBc + V = L$ to ensure that $\diamondsuit$ holds and that we have ground model definability. Take a generic $H \subseteq \Add(\omega_1, 1)$ and consider $W^V_H$ from \myref{th:DC_deltaHOldsDC_delta++Fails}. This is a model of $\GBcminus + \DC_{\omega_1}$-scheme in which $\mathcal{P}(\omega_1)$ does not surject onto $\omega_3$. By \myref{prop:chainConditionPretame}, since the $\DC_{\omega_1}$-scheme holds, any ccc class forcing over $W^V_H$ is pretame. Specifically, the class tree iteration $\p(\vec \J,\Ord^{{<}\omega})$ remains pretame in $W^V_H$. Next, note that $H$ added no new subsets of $\omega$, which in particular means that $\diamondsuit$ holds in $W^V_H$. Thus $\p(\vec \J,\Ord^{{<}\omega})$ satisfies all the necessary properties mentioned in \myref{subsec:Jensen} and, for any generic $G$, $W^V_H$ is definble in $W^V_H[G]$. So, let $G \subseteq \p(\vec \J,\Ord^{{<}\omega})$ be generic and consider $W = \bigcup_{\mathcal{T} \in \mathbb{T}} W^V_H[G_\mathcal{T}]$ as before. By the previous analysis, it is clear that $W \models \ZFCminus + \neg \DC_\omega$-scheme.

Finally, while $\mathcal{P}(\omega_1)$ is now a big proper class in $W$, $\mathcal{P}(\omega_1) \cap W^V_H$ is not. This is because $\mathcal{P}(\omega_1) \cap W^V_H$ is a set of cardinality $\omega_2$ in $L[H]$ and, since the second forcing doesn't collapse cardinals, this must still be true in $L[H][G]$. Therefore this class cannot possibly surject onto $\omega_3$ in $W$. Thus, $W$ is our desired model of $\ZFCminus$ with a proper class that is not big.
\end{proof}

\section{Embeddings with \texorpdfstring{$\mathcal P (\omega)$}{P(omega)} a proper class}

\noindent In this section, we show that there is a model of $\ZFCminus$ having a definable elementary embedding $j$ with a critical point in which $\mathcal P (\omega)$ does not exist.

Suppose that $V\models\ZFC$ and $\kappa$ is a measurable cardinal in $V$ with a normal measure $U$. Let $j \colon V\to M$ be the ultrapower map by $U$. Let $\p=\Add(\omega,1)$, $\q=\suppprod{\p}{\omega}{\omega}$, and $G\subseteq \q$ be $V$-generic. We construct $W^V_G$ in $V[G]$ as in \myref{sec:Zarach}.

A folklore result, known as the Lifting Criterion, states that given an elementary embedding $j \colon V\to M$, a poset $\p\in V$, a $V$-generic filter $G\subseteq \p$ and an $M$-generic filter $H\subseteq j(\p)$, we can lift (extend) the embedding $j$ to $j \colon V[G]\to M[H]$ if and only if $j\image G\subseteq H$. In case we can lift, the lift $j$ is given by $j(\dot x_G)=j(\dot x)_H$. Thus, by the Lifting Criterion, using that $j(\q)=\q$ and $j\image G=G$, we can lift $j$ to the embedding $j \colon V[G]\to M[G]$. Moreover, we can show that $j$ is the ultrapower map by the measure $U_\omega$ generated by $U$, namely $A\in U_\omega$ if and only if there is $B\in U$ such that $B\subseteq A$. This also shows that, for every $n<\omega$, $j$ lifts to $j_n \colon \unionmodel{V}{n}\to \unionmodel{M}{n}$ ($\unionmodel{V}{n}=V[G_n]$, $\unionmodel{M}{n}=M[G_n]$), and is the ultrapower map by the measure $U_n$ generated by $U$ in $\unionmodel{V}{n}$. Thus,
\[
U\subseteq U_1\subseteq U_2\subseteq \cdots\subseteq U_n\subseteq \cdots\subseteq U_\omega.
\]

By the definability of $W^V_G$ in $V[G]$, we can restrict the lift $j \colon V[G]\to M[G]$ to an elementary embedding
\[
j^W \colon W^V_G\to W^M_G,
\]
where $W^M_G$ is constructed in $M[G]$ analogously to $W^V_G$.

\begin{proposition}
$j^W=\bigcup_{n<\omega}j_n$.
\end{proposition}

\begin{proof}
Observe that it suffices to show that the lift $j_n \colon \unionmodel{V}{n}[G_{n,\tail}]\to \unionmodel{M}{n}[G_{n,\tail}]$ is the lift $j \colon V[G]\to M[G]$. But this is clear because, by our observation above, the lift of $j_n$ is the ultrapower map by the measure generated by $U_n$, which is clearly $U_\omega$.
\end{proof}

\noindent Let ${\rm E}$ be the membership relation modulo $U_\omega$. Observe that $U_\omega$-equivalence and ${\rm E}$ are both definable in $W^V_G$ from the set $U\in W^V_G$. Fix a function $f \colon \kappa\to W^V_G$ in some $\unionmodel{V}{n}$.  Since $j^W=\bigcup_{n<\omega}j_n$, we know that whenever $g \colon \kappa\to W^V_G$ from $W^V_G$ is an ${\rm E}$-member of $f$, then $g$ is $U_\omega$-equivalent to some $g^* \colon \kappa\to \unionmodel{V}{n}$ with $g^*\in \unionmodel{V}{n}$. Thus, $\unionmodel{V}{n}$ has a set $X^f$ consisting of functions $g \colon \kappa\to \unionmodel{V}{n}$ that are ${\rm E}$-members of $f$ such that any function $h \colon \kappa\to W^V_G$ from $W^V_G$ that is an ${\rm E}$-member of $f$ is $U_\omega$-equivalent to a function in $X^f$. Thus, in $W^V_G$, given any function $f \colon \kappa\to W^V_G$, we can associate to it the set $X^f$.

Now let's provide a definition of $j$ in $W^V_G$. Fix a function $f \colon \kappa\to W^V_G$, and let $b\in V[G]$ be the image of $[f]_{U_\omega}$ under the transitive collapse. Working in $W^V_G$, let $X^f_0=X^f$. Now suppose inductively that we are given $X^f_n$, and let $X^f_{n+1}=\bigcup_{g\in X^f_n}X^g$. Let $X^f_\omega=\bigcup_{n<\omega}X^f_n$. It should be clear that the transitive collapse of $\la X^f_\omega,{\rm E}\ra$ (modulo $U_\omega$) is the transitive closure of $b$, from which we can compute $b$. Thus, $W^V_G$ can compute $j(a)$ by computing $X^{c_a}_\omega$, where $c_a \colon \kappa\to \{a\}$ is the constant function.

Finally, recall that crucial property that led to the failure of the $\DC_{\omega_1}$-scheme in \Cref{th:DC_omegaHoldsDC_omega_2Fails} was that the ground model ($L$) was definable in $W^V_G$. Therefore, if $V = L[U]$ then we will again have that the $\DC_{\omega_1}$-scheme fails in the resulting model.

Putting together all of the above, we obtain the following result.

\begin{theorem}\label{th:EmbeddingWithoutPowerSet}\phantom{a}
\begin{enumerate}
\item $W^V_G\models\ZFCminus+\DC_\omega$-scheme.
\item The $\DC_{\omega_2}$-scheme fails in $W^L_G$.
\item $\mathcal P (\omega)$ (and therefore $V_\alpha$ for $\alpha > \omega$) does not exist in $W^V_G$.
\item $W^V_G$ has a definable elementary embedding $j$ with a critical point.
\item If additionally $V = L[U]$, then the $\DC_{\omega_1}$-scheme fails in $W^L_G$.
\end{enumerate}
\end{theorem}

\noindent We can also use the construction of \myref{sec:Zarach} to produce other examples of the limitations to \myref{th:elementaryEmbeddingHierarchy}. For instance, it is easy to produce a model $W\models\ZFCminus$ with a cofinal elementary embedding $j \colon W\to W$ having a critical point (but $W$ won't satisfy $\ZFCminus_j$).

Start with a transitive model $M\models\ZFC$ for which there is an elementary embedding $j \colon M\to M$ with critical point some ordinal $\kappa$ (the consistency strength of this assumption is below $0^{\#}$). Let $\p=\Add(\omega,1)$ in $M$, force with $\q=\suppprod{\p}{\omega}{\omega}$, and let $G\subseteq \q$ be $M$-generic. First, we lift $j$ to an elementary embedding $j \colon M[G]\to M[G]$, and then restrict $j$ to $j^W \colon W^M_G\to W^M_G$.

We can also use the construction of \myref{sec:Zarach} to produce a model $W\models\ZFCminus_j$ such that $\mathcal{P}(\omega)$ is a proper class and where $j \colon W\to W$ is an elementary embedding with a critical point, but not cofinal.

Suppose that $V\models\ZFC+\Ione$ and fix an elementary embedding
\[
j \colon V_{\lambda+1}\to V_{\lambda+1}
\]
with critical point $\kappa<\lambda$. Let
\[
\jextended \colon H_{\lambda^+}\to H_{\lambda^+}
\]
be the elementary embedding obtained from $j$. Let $\p=\Add(\omega,1)$, $G\subseteq \q=\suppprod{\p}{\omega}{\omega}$ be $V$-generic, and $W^V_G$ be constructed in $V[G]$ as in \myref{sec:Zarach}. Let
\vspace{-10pt}
\[
\jextended \colon H_{\lambda^+}[G]\to H_{\lambda^+}[G]
\]
be the lift of $\jextended$ given by $\jextended(\dot x_G)=\jextended(\dot x)_G$. Let
\[
N^V_G=\bigcup_{n<\omega}H_{\lambda^+}[G_n]
\]
be the proper class of $W^V_G$ consisting of sets whose transitive closure has size at most $\lambda$. Let
\[
j^N \colon N^V_G\to N^V_G
\]
be the restriction of $\jextended$ to $N^V_G$. Then the embedding $j^N$ can be defined by the formula $\varphi^j(x,y,\p,j,\lambda,G)$ asserting that:
\begin{itemize}
\item $x\in W^V_G$ has transitive closure of size at most $\lambda$,
\item there is a $\q$-name $\dot x\in H_{\lambda^+}^V$ such that $x=\dot x_G$, and
\item $y=j^+(\dot x)_G$.
\end{itemize}
Next, let's argue that if $\pi$ is any coordinate-switching automorphism of $\q$, then $\varphi^j(x,y,\p,j,\lambda,\pi\image G)$ also defines $j^N$. Suppose that for some $a\in N^V_G$, $\jextended(a)=b$. Let $a=\dot x_G$ with $\dot x\in H_{\lambda^+}^V$. Then $b=\jextended(\dot x)_G$. Clearly, $a= (\pi (\dot x))_{\pi\image G}$. Now,
\[
\jextended(\pi(\dot x)) = \jextended(\pi)(\jextended(\dot x)) = \pi(\jextended(\dot x)).
\]
Thus, $\jextended(\pi (\dot x))_{\pi\image G}=\jextended(\dot x)_G=b$. Thus, using $G$ or $\pi\image G$ both yield the same embedding $j^N$. From this we immediately get the following strengthening of the key elementarity lemma from \myref{sec:Zarach}, where we let $\jextended_n \colon H_{\lambda^+}^{\unionmodel{V}{n}}\to H_{\lambda^+}^{\unionmodel{V}{n}}$ be the lift of $\jextended$ in $\unionmodel{V}{n}$ and $j^N_n$ be its restriction to $N^{\unionmodel{V}{n}}_{G^{(n)}}$.

\begin{lemma}\label{lem:elementarityWwithEmbedding}
For every $n<\omega$, $\la W^{\unionmodel{V}{n}}_{G^{(n)}},j_n^{N}\ra\prec \la W^V_G,j^N\ra$.
\end{lemma}

\noindent Using the arguments from \myref{sec:Zarach}, \myref{lem:elementarityWwithEmbedding} gives:

\begin{theorem}
$W^V_G\models \ZFCminus_{j^N}$, and hence $N^V_G \models \ZFCminus_{j^N}$.
\end{theorem}

\begin{proof}
It follows from \myref{lem:elementarityWwithEmbedding} that $W^V_G\models \ZFCminus_{j^N}$. It remains to observe that any instance of collection for a set whose transitive closure has size at most $\lambda$ can also be assumed to have transitive closure of size at most $\lambda$.
\end{proof}

\subsection*{Acknowledgements} An early version of some of the ideas in this paper originally appeared as part of the second author's PhD thesis, supervised by Michael Rathjen and Andrew Brooke-Taylor. The second author was supported by the UK Engineering and Physical Sciences Research Council during the research for this work and is grateful for their support.

\bibliography{UnionModels}
\bibliographystyle{alpha}

\end{document}